\documentclass[a4paper,11pt]{amsart}

\usepackage{ucs}
\usepackage[latin9]{inputenc}
\usepackage{graphicx}
\usepackage{amsfonts}
\usepackage{dsfont}
\usepackage{amssymb}
\usepackage{amsmath}
\usepackage{amsthm}
\usepackage{enumerate}
\usepackage{stmaryrd}
\usepackage{fullpage}
\usepackage{ifthen}
\usepackage{subfigure}
\usepackage{epic}
\usepackage{textcomp}
\usepackage{mathrsfs}

\usepackage[hypertexnames=false,colorlinks=true,linkcolor=blue,citecolor=blue]{hyperref}
\usepackage[numbers,comma,square,sort&compress]{natbib}

\usepackage{color}

\graphicspath{{eps/}{pdf/}}

\newcommand{\bqq}{\begin{equation}}
\newcommand{\eqq}{\end{equation}}
\newcommand{\bqs}{\begin{equation*}}
\newcommand{\eqs}{\end{equation*}}

\DeclareMathOperator{\sgn}{sign} 

\newcommand{\A}{\mathbb{A}}
\newcommand{\B}{\mathbb{B}}
\newcommand{\D}{\mathbb{D}}

\newcommand{\R}{\mathbb{R}} 
\newcommand{\C}{\mathbb{C}}
\newcommand{\N}{\mathbb{N}}
\newcommand{\G}{\mathbb{G}}
\newcommand{\MM}{\mathbb{M}}
\newcommand{\LL}{\mathbb{L}}

\newcommand{\md}{\mathrm{d}}
\newcommand{\mbi}{\mathbf{i}}

\newcommand{\bQ}{\mathbb{Q}}
\newcommand{\bI}{\mathbb{I}}
\newcommand{\cQ}{\mathscr{Q}}

\newcommand{\K}{\mathcal{K}}

\renewcommand{\P}{\mathcal{P}}

\renewcommand{\hom}{\mathrm{hom}}
\newcommand{\bc}{\mathrm{bc}}
\newcommand{\dir}{\mathrm{dir}}
\newcommand{\rfl}{\mathrm{ref}}
\newcommand{\trp}{\mathrm{trap}}

\newcommand{\Ran}{\mathrm{Ran}}

\newcommand{\M}{\mathcal{M}}
\newcommand{\bU}{\underline{U}}

\renewcommand{\Re}{\mathrm{Re}}

\newtheorem{lem}{Lemma}[section]
\newtheorem{thm}{Theorem}
\newtheorem{prop}[lem]{Proposition}

\newtheorem{rmk}[lem]{Remark}

 {\begin{trivlist}\item[]\textbf{Proof#1 }}%
 {\hspace*{\fill}$\rule{0.3\baselineskip}{0.35\baselineskip}$\end{trivlist}}

  {\begin{trivlist}\item[]{\bf Hypothesis #1 }\em}{\end{trivlist}}

\numberwithin{equation}{section}

\title[Riemann shocks in systems]{Exponential asymptotic stability of Riemann shocks\\of hyperbolic systems of balance laws}

\author{Gr\'egory Faye}
\address{Institut de Math\'ematiques de Toulouse ; UMR5219, Universit\'e de Toulouse ; CNRS UPS IMT, F-31062 Toulouse Cedex 9 France}
\email{{\tt gregory.faye@math.univ-toulouse.fr}}
\thanks{G.F. acknowledges support from the ANR via the project Indyana under grant agreement ANR- 21- CE40-0008, Labex CIMI under grant agreement ANR-11-LABX-0040, and an ANITI (Artificial and Natural Intelligence Toulouse Institute) Research Chair.}

\author{L.~Miguel Rodrigues}
\address{
Univ Rennes \& IUF, CNRS, IRMAR - UMR 6625, F-35000 Rennes, France}
\email{{\tt luis-miguel.rodrigues@univ-rennes1.fr}}
\thanks{}

\begin{document}

\maketitle

\begin{abstract}
\noindent For strictly entropic Riemann shock solutions of strictly hyperbolic systems of balance laws, we prove that exponential spectral stability implies large-time asymptotic orbital stability. As a preparation, we also prove similar results for constant solutions of initial value and initial boundary value problems, that seem to be new in this generality. Main key technical ingredients include the design of a nonlinear change of variables providing a hypocoercive Kawashima-type structure with dissipative boundary conditions in the high-frequency regime and the explicit identification of most singular parts of the linearized evolution, both being deduced from the mere spectral assumption. 
\end{abstract}


 

\section{Introduction}

\subsection{Overall motivation}

The present contribution brings a major piece to the still ongoing development of a Lyapunov theory for traveling waves of hyperbolic systems. By this we mean a theory that describes the large-time dynamics near spectrally stable waves in arbitrarily strong topologies.

Surprisingly enough, though modeling with hyperbolic systems is almost ubiquitous and the understanding of nonlinear waves is an important part of the qualitative analysis of any system, the hyperbolic nonlinear wave stability theory is still in its infancy by many respects. In particular, with the notable exception of waves of scalar balance laws \cite{DR1,DR2,GR} and discontinuous fronts of a specific $2\times 2$ system considered in \cite{Yang-Zumbrun}, analyses in the literature fail to cover nonlinear waves with discontinuous\footnote{With this respect, we warn the reader against the confusion that could arise from the fact that on one hand, there is a large body of literature using the terminology \emph{shock} to name some smooth fronts, and on the other hand, there is an equally large body of literature studying shocks, but using the word \emph{stability} in the sense of \emph{short-time} persistence.} piecewise-smooth profiles, or, even, those with smooth profiles but a characteristic point. It is all the more regrettable that the emergence of such kind of objects is a distinctive feature of hyperbolic systems, and sometimes the main reason to adopt an hyperbolic model rather than a parabolic one.

The main explanation for this gap in the theory is obviously that classical stability theory for one-dimensional traveling waves--- as described for instance in \cite{Kapitula-Promislow,JNRZ} ---, or more generally for radial or planar multidimensional waves, use both the regularity of wave profiles and the elliptic nature of non-characteristic one-dimensional operators at almost every stage of the analysis. Consistently, in the reverse direction, we point out that for waves of hyperbolic systems with non-characteristic smooth profiles, a rather comprehensive theory is indeed available; see for instance \cite{Mascia-Zumbrun-1,Mascia-Zumbrun-2}.

Up to our knowledge, the present contribution is the first one to provide a Lyapunov-type stability result for some discontinuous waves of a large class of hyperbolic systems, including systems of arbitrary dimension. With this respect, it is important to note that the structure of scalar equations or $2\times 2$ systems is highly non representative of the general structure of hyperbolic systems; see detailed discussions in \cite{Serre-conservation-I,Serre-conservation-II,Bressan,Benzoni-Serre} or \cite[Appendix~A3]{MR}. A related observation is that spectral stability of waves considered in \cite{DR1,DR2,GR,Yang-Zumbrun} is fully elucidated\footnote{In the sense that it is reduced to sign conditions on a few numbers.} either directly along the linear stability analysis or in a companion paper (\cite{SYZ} for \cite{Yang-Zumbrun}) instead of being taken as an abstract assumption. The cases when this is possible analytically are quite exceptional, even at the ODE level. Roughly speaking, scalar equations and $2\times 2$ systems are for hyperbolic systems as exceptional as scalar reaction-diffusion equations --- analyzable through Sturm-Liouville theory and maximum principles --- are with respect to general parabolic systems.

However, in order to contain technicalities as much as possible in this first contribution to the general system theory, we do make a few restrictions in generality. First, we restrict to strictly hyperbolic systems. We expect that though it would be interesting to relax these assumptions so as to enlarge the range of applicability of the results, this extension would not bring any new dynamical phenomenon. 

The most stringent restrictions we consider are on the class of traveling waves we study. Firstly, we focus on profiles that are piece-wise constant with a single discontinuity, so-called Riemann shocks. Secondly, we assume that involved profiles are non characteristic. At last, we only tackle the case when spectral stability holds with a spectral gap, hence is expected to yield time-exponential decay (in a suitable orbital sense), as opposed to algebraic decay.  

Incidentally we point out that, as we show below, our restrictions imply that the discontinuities of profiles we consider are of strictly entropic Lax type but not necessarily of extreme Lax type\footnote{As long as it is not essential to our analysis, in order to maintain reading fluidity we shall not define explicitly neither standard hyperbolic terminology --- such as Lax shock --- nor standard wave terminology. We refer the reader to \cite{Serre-conservation-I,Serre-conservation-II,Bressan,Benzoni-Serre} on the former and to \cite{Kapitula-Promislow} on the latter.}.

\subsection{Main statement}

To be more concrete, we consider a system of balance laws
\bqq
\partial_t U+ \partial_x(A(U))=g(U)
\label{edp_intro}
\eqq
with smooth coefficients $A$, $g$, and unknown $U$ depending on time variable $t\in\R$, space variable $x\in\R$ and taking values in $\R^n$, $n\in\N$.

We are interested in the dynamics near a traveling-wave solution $(t,x)\mapsto \bU(x-\sigma t)$ with speed $\sigma\in\R$ and profile $\bU$, of Riemann-shock type
\bqs
\bU(\xi)=\left\{
\begin{array}{ll}
\bU_-, & \text{ if } \xi<0, \\
\bU_+,& \text{ if } \xi>0,
\end{array}
\right.
\eqs
for some $\bU_+\in\R^n$, $\bU_-\in\R^n$. The fact that this is indeed a weak solution to \eqref{edp_intro} is equivalent to $\bU_+$, $\bU_-$ being equilibria, that is,  $g(\bU_\pm)=0_{\,\R^n}$, and $\bU_+$, $\bU_-$ being connected by Rakine-Hugoniot condition at speed $\sigma$
\bqq\label{eq:RH_intro}
A(\bU_+)-A(\bU_-)=\sigma \left(\bU_+-\bU_- \right).
\eqq

\begin{rmk}\label{rk:count}
The reader may rightfully wonder whether the object under study is structurally stable, or in more quantitative words, may ask how many parameters should the system contain to guarantee that if one perturbs the system then a similar object exists for some nearby parameters. Our spectral assumptions, to be detailed below, imply that this count is the same as the count of solutions to \eqref{eq:RH_intro}, when $\bU_+$ and $\bU_-$ are fixed (here determined\footnote{This implicitly takes for granted that involved zeros of $g$ are non-degenerate. But this is indeed a consequence of exponential spectral stability assumed below.} as zeros of $g$) but $\sigma$ is free. Therefore, with $n$ equations, $(n-1)$ parameters are needed to ensure structural persistence. We point out that for similarly exponentially stable waves with a single shock, that are only constant on one side (respectively constant on no side), a similar count would indicate that $(n-2)$ parameters (resp. $(n-3)$ parameters) are needed. We have left for a further contribution the study of such objects, because we expect that, though their stability may be analyzed with similar techniques, this would introduce an extra layer of complexity in an already quite technical proof.  
\end{rmk}

As is well-known, one should impose some extra conditions, of entropy type, to single out unique solutions among the otherwise large set of weak solutions. Our spectral stability assumptions do imply that $(t,x)\mapsto \bU(x-\sigma t)$ satisfies strict Lax entropy conditions. But even this, alone, is not sufficient to conclude uniqueness from known results. Indeed, further assumptions both on the structure of the system (strict hyperbolicity, genuine nonlinearity or linear degeneracy of characteristic fields,...) and/or on the solutions (small BV norm, Riemann data, piecewise smoothness,...) are involved in the classical uniqueness theory. On the latter\footnote{We warn the reader that many parts of the literature restrict to conservations laws, for the sake of simplicity of exposition. However, as far as short-time existence is concerned, the extension to balance laws is straightforward.}, we refer to \cite{Li-Yu,Bressan} for results specific to one-dimensional solutions (as considered here) and \cite{Majda1,Majda2,Metivier,Benzoni-Serre} for (partial) multidimensional counterparts. In the present contribution, we only use results for piecewise smooth solutions of strictly hyperbolic systems. This choice implicitly hinges on the expectation that whatever choice is made to ensure uniqueness, if the criterion holds in a strict sense for the background wave it will still hold for nearby solutions built with the same regularity structure ;  see the explicit scalar discussion in \cite{DR1,DR2}.

We assume that the system is strictly hyperbolic and non-characteristic, near $\bU$. This amounts to enforcing that both $D_UA(\bU_-)$ and $D_UA(\bU_+)$ have $n$ distinct real eigenvalues, and none of those are equal to $\sigma$. We expect that the strict hyperbolicity assumption could be relaxed in a relatively standard way, but to the price of cumbersome extra technical details.

To state what spectral stability means, we first generalize the discussion about what this means to be a weak solution to more general functions, with the same structure of regularity. Consider a locally bounded function $U$ defined by
\bqs
U(t,x)=\bU(x-\sigma t -\psi(t))+V(t,x-\sigma t -\psi(t)),
\eqs
with say $V$ $\mathscr{C}^1$ on $\R\times \R^*$ with limits from the right and from the left on $\R\times\{0\}$ and $\psi$ $\mathscr{C}^1$, encoding respectively perturbations in shape and in position of the Riemann shock. Such a $U$ solves weakly \eqref{edp_intro} if and only if 
\bqq
\left\{
\begin{split}
\partial_t V+ \left(D_UA(\bU+V)-(\sigma +\psi'(t))I_n\right)\partial_x V&=g(\bU+V),\qquad (t,x)\in\R\times \R^*, \\
\quad-\left(\sigma+\psi'(t)\right)\left[\bU+V(t,\cdot)\right]_{0}+\left[A(\bU+V(t,\cdot))\right]_{0}&=0\,,\qquad\qquad\qquad t\in\R\,.
\end{split}
\right.
\label{systedp_intro}
\eqq
In the latter we have used jump notation  $[W]_0:=W(0^+)-W(0^-)$, where $W(0^\pm)=\lim_{h\searrow0}W(0\pm h)$. 

Linearizing System~\eqref{systedp_intro} in $(V,\psi)$ small leaves 
\bqq
\left\{
\begin{split}
\partial_t V&+ \A_\pm \partial_x V=\G_\pm V,\qquad\textrm{on }\R\times \R_\pm, \\
-\psi'\left[\bU\right]_{0}&+\left[(D_UA(\bU)-\sigma I_n)V\right]_{0}=0,
\end{split}
\right.
\label{linsystedp_intro}
\eqq
where 
\begin{align*}
\A_\pm&:=D_UA(\bU_\pm)-\sigma I_n\,,&
\G_\pm&:=D_U G(\bU_\pm)\,.
\end{align*}
For any spectral parameter $\lambda\in\C$, this yields as a spectral problem
\bqq
\left\{
\begin{split}
\left(\lambda I_n+\A_\pm \partial_x-\G_\pm\right)&\widetilde{V}_\pm=F_\pm, \qquad 
\textrm{on }\R_\pm,\\
-\lambda \widetilde{\psi}\left[\bU\right]_{0}+\A_+\widetilde{V}_+(0)&-\A_-\widetilde{V}_-(0)=F_0.
\end{split}
\right.
\label{linspectral_intro}
\eqq
Choosing $L^2$ as a reference functional space, we say that $\lambda$ does not belong to the spectrum of \eqref{linsystedp_intro} if and only if for any $F_\pm\in L^2(\R_\pm;\C^n)$ and any $F_0\in\C^n$, there exists a unique $(\widetilde{V}_+,\widetilde{V}_-,\widetilde{\psi})\in H^1(\R_+;\C^n)\times H^1(\R_-;\C^n)\times \C$ solving \eqref{linspectral_intro}. Note that by the uniform boundedness principle, in this case, $(\widetilde{V}_+,\widetilde{V}_-,\widetilde{\psi})$ depends boundedly on $(F_+,F_-,F_0)$, so that solving \eqref{linspectral_intro} defines a bounded linear operator, called resolvent operator. One readily checks that the spectrum forms a closed set and that on the complementary of the spectrum, the resolvent map depends analytically on $\lambda$. Consistently, we shall say that an element $\lambda_0$ of the spectrum has finite multiplicity if $\lambda_0$ is isolated in the spectrum, and the resolvent map possesses a meromorphic singularity at $\lambda_0$, with finite-rank residue, the algebraic multiplicity of $\lambda_0$ being then defined as this rank.

The problem is invariant by spatial translation, that is, any translate of a solution is still a solution. Since $\bU$ is not constant, thus not invariant by spatial translations, this causes $0$ to be in the spectrum. Explicitly, when $\lambda=0$, $(\widetilde{V}_+,\widetilde{V}_-,\widetilde{\psi})\equiv(0_{\,\R^n},0_{\,\R^n},1)$ defines a non-zero solution to \eqref{linspectral_intro} with $(F_+,F_-,F_0)\equiv(0_{\,\R^n},0_{\,\R^n},0)$. At the spectral level, the best one can expect is therefore that the only obstacle to exponential decay is a simple eigenvalue at $0$. 

At the nonlinear level, one could expect then, consistently with classical analysis for smooth waves, that under this spectral stability condition, exponential orbital stability holds, that is exponential decay of the distance to the family of spatial translates of $\bU$ does occur. However, this can happen only for initial data compatible with the regularity structure of $\bU$, not only in the sense that they are piecewise smooth with a single discontinuity but also in the sense that they are compatible with the short-time persistence of a single-shock structure. In other words, we prove global-in-time stability only for perturbations for which short-time stability holds (in the sense of \cite{Majda2,Metivier}). By a dimensional count similar to the one carried out in Remark~\ref{rk:count}, one may check that this is (at least\footnote{The exact number depends on the level of regularity one wants to enforce.}) an $(n-1)$-dimensional constraint on the discontinuity, that would otherwise be resolved by generating the superposition of $n$ simple waves (at least in simplest cases). A simple way, used in \cite{Yang-Zumbrun}, to impose this constraint is to restrict to perturbations that are supported away from the discontinuity. Beyond simplicity, the foregoing way also comes with the advantage that such data are prepared to propagate any level of regularity.

The content of our main theorem is precisely to prove exponential orbital stability under small perturbations that do not disintegrate the shock instantaneously, provided that the wave under consideration is exponentially spectrally stable. To keep the statement as streamlined as possible, we introduce beforehand the terminology that an initial shape perturbation $V_0$ is $H^2$-compatible if there exist $(\psi_1,\psi_2)\in\R^2$ such that $\left[A(\bU+V_0)\right]_{0}=\left(\sigma+\psi_1\right)\left[\bU+V_0\right]_{0}$ and
\begin{align*}
\left[\left(\,D_UA(\bU+V_0)-\left(\sigma+\psi_1\right)\,I_n\,\right)
\left(-\partial_x(A(\bU+V_0)-\left(\sigma+\psi_1\right)V_0)+g(\bU+V_0)\right)\right]_{0}&=\left(\sigma+\psi_2\right)\left[\bU+V_0\right]_{0}\,.
\end{align*}

\begin{thm}\label{th:main}
Assume that there exists $\alpha_0>0$ such that the spectrum of \eqref{linsystedp_intro} is contained in 
\[
\{\,\lambda\in\C\,;\,\Re(\lambda)<-\alpha_0\,\}\,\cup\,\{0\}
\]
and that $0$ is a simple eigenvalue. For any $0<\alpha<\alpha_0$, there exist positive $C_0$ and $\epsilon_0$ such that for any $(V_0,\psi_0)\in H^2(\R^*;\R^n)\times \R$ with $\|V_0\|_{H^2(\R^*)}\leq \epsilon_0$ and $V_0$ $H^2$-compatible, there exist $V\in\mathscr{C}^0(\R_+;H^2(\R^*;\R^n))\cap \mathscr{C}^1(\R_+;H^1(\R^*;\R^n))$ and $\psi\in \mathscr{C}^2(\R_+)$ such that 
\[
(t,x)\mapsto \bU(x-\sigma t -\psi(t))+V(t,x-\sigma t -\psi(t))
\]
solves \eqref{edp_intro} with initial data $(\bU+V_0)(\cdot-\psi_0)$ and
\begin{align*}
\|V(t,\cdot)\|_{H^2(\R^*)}
\,+\,|\psi'(t)|\,+\,|\psi''(t)|
&\leq\,C_0\,e^{-\alpha\,t}\,\|V_0\|_{H^2(\R^*)}\,,&t\in\R_+\,,
\end{align*}
and for some $\psi_\infty\in\R$,
\begin{align*}
|\psi_\infty-\psi_0|
&\leq\,C_0\,\|V_0\|_{H^2(\R^*)}\,,&\\
|\psi(t)-\psi_\infty|
&\leq\,C_0\,e^{-\alpha\,t}\,\|V_0\|_{H^2(\R^*)}\,,&t\in\R_+\,.
\end{align*}
\end{thm}

The regularity threshold for piecewise-smooth solutions is piecewise Lipschitz regularity. We have decided to work with $L^2$-based spaces and integer-valued derivatives, hence the choice of an $H^2$ framework. However our proof would also provide similar results in $W^{s,p}$ spaces, provided that $1\leq p<\infty$, $s>1+1/p$. In contrast, we expect that it cannot be easily adapted to reach the sharp $W^{1,\infty}$ stage, at which the scalar analysis is performed in \cite{DR1,DR2,GR,BR}. It is also worth pointing out that we could also prove exponential decay for norms encoding any extra regularity assumed on $V_0$ without imposing any extra smallness condition, as in \cite[Propagation~3.5]{DR1}. We mention however that the notion of compatibility of $V_0$ with the single-shock structure should be adapted accordingly, the number of constraints increasing with the degree of regularity. 

The spectral stability assumption is formulated in terms of spectral problems to appear as natural as possible and to highlight that it is sharp at the linearized level, as a Lyapunov-type result should be. However in the very first step of our proof we show that it implies that our background shock is of Lax type and that constant equilibria $U\equiv \bU_+$ and $U\equiv \bU_+$ are themselves exponentially spectrally stable, and, assuming the latter, that it is equivalent to a non-vanishing condition on an Evans-Lopatinski\u{\i} determinant, a form more commonly encountered in the literature about spectral and linear stability of shocks (see for instance \cite{Godillon,Godillon-Lorin,Texier-Zumbrun,JNRYZ}). 

Concerning spectral stability assumption, we stress that it is not designed to be easily checked analytically but to be sharp. More, in the few cases where one could prove spectral stability (dissipative rich systems, small-amplitude waves,...), it is also reasonable to expect, and common in related situations, that the spectral argument could be upgraded\footnote{Quite often only the nonlinear argument appears in the literature.} into a direct simpler proof of nonlinear stability. See for instance \cite[Appendix~A]{R_HDR} and \cite{RZ} for an explicitly worked-out correspondence, restricted to high-frequency stability though. However, from a more applied point of view, the reduction to an Evans-Lopatinski\u{\i} determinant condition brings the spectral stability issue to a stage reasonably decidable by well-conditioned numerics. See the related periodic wave study in \cite[Section~7]{JNRYZ}, that expands on algorithms initially developed for smooth waves (see \cite[Appendix~D]{BJNRZ} and \cite[Chapter~3]{R_HDR}).

For the convenience of the reader, let us sketch the main features of the proof of Theorem~\ref{th:main}. It relies on two kinds of nonlinear estimates. On one hand, we prove that spectral stability implies the existence of a nonlinear change of coordinates adapted to the high-frequency regime, yielding an hypocoercive Kawashima-type structure with dissipative boundary conditions for higher-order derivatives, with forcing terms due to small-order derivatives. This results in nonlinear estimates proving that the decay of higher-order derivatives is slaved to the one of low-order ones. On the other hand, we estimate low-order derivatives \emph{via} Duhamel formula. This causes an apparent loss of derivatives that is cured with the above-mentioned high-frequency damping estimates. The key to these low-regularity estimates is a careful study of the linearized dynamics. Concerning the latter, we mention that our analysis of the inversion of Laplace transforms, required to go from spectral problems to linearized time-evolution, relies on the explicit computation of the singular parts of the dynamics up to an order where the inversion becomes regular. A fine description of high-frequency expansions of the spectral problems are obviously used here crucially and explicitly, but those are actually also a key-point of the nonlinear high-frequency damping estimates.

\begin{rmk}\label{rk:energy}
It is part of the standard Lyapunov theory that for constant equilibria of finite-dimensional differential equations, exponential spectral stability implies the existence of a nonlinear energy estimate, sufficient to deduce exponential nonlinear stability. Our proof extends such a philosophy to the high-frequency regime of our current problem. Some readers may rightfully wonder whether, likewise, a full stability result could be obtained by a pure nonlinear energy estimate. An expected gain is that proofs by pure energy estimates tend to be less technical than our two-tiers proof. However, as we show in Appendix~\ref{s:example}, in general this expectation cannot be met even for constant equilibria of \eqref{edp_intro}, except for scalar equations and systems of two equations. In other words, even in the simplest infinite-dimensional cases, the nonlinear stability results deduced from spectral stability assumptions cover more cases than those that may be proved by energy estimates, the price to pay being technical complexity of the required proof. 
\end{rmk}

\subsection{Outline and perspectives}

The present contribution focuses on the simplest non trivial system case. Yet, for applications, there are a quite large series of extensions that are worth carrying out, including the consideration of profiles that are not piecewise constant but still asymptotically constant, of periodic profiles, of profiles with characteristic points, of stability as plane waves of multidimensional systems, of cases when zeros of the source term $g$ (thus constant equilibria) are not isolated but form instead a smooth manifold,...

Let us give only a few hints about difficulties and novelties to be expected from these desired extensions. Concerning characteristic points, we stress that the scalar analysis is already contained in \cite{DR2} and it reveals a dramatic influence on the nature of the spectral problem\footnote{For instance, in the characteristic case, the spectrum depends crucially on the chosen level of regularity encoded in the underlying functional space.} and on the phase dynamics. As for multidimensional plane Riemann shock stability, already in the scalar case dealt with in \cite{DR1} one derives that the spectrum necessarily\footnote{Or the eigenvalue $0$ has infinite-dimensional multiplicity.} includes the whole imaginary axis, and correspondingly perturbations on the shape of the shock location do not flatten back. At last, we point out that both the periodic case and the case when $g$ is not full-rank preclude any spectral gap, so that only algebraic decay is to be expected.

The rest of the present paper is devoted to the proof of Theorem~\ref{th:main}. However, for expository reasons, instead of focusing directly on it, we consider simpler problems so as to gradually introduce technical arguments and conclude with the proof of Theorem~\ref{th:main}. Explicitly, we consider as intermediate steps towards our main goal the stability of constant equilibria first for the initial value problem associated with equations posed on the whole line then for initial boundary value problems posed on half-lines. Incidentally let us observe that a similar choice was made in \cite{Metivier_Bourbaki} to expound the content of \cite{Majda1,Majda2}. Though we provide these other stability results as intermediate expository steps, they seem to be new in this generality.

\medskip

\noindent \emph{Acknowledgments:} M.R. expresses his gratitude to Vincent Duch\^ene for enlightening discussions at an early stage of the project. G.F. thanks IRMAR and M.R. thanks I.M.T. for their hospitality during respective visits.

\section{Stability of constant equilibria}\label{s:constant}

We begin by revisiting the Lyapunov stability theory for constant solutions. Unlike the analysis of the scalar case in \cite{DR1}, we shall not use the corresponding result to prove Theorem~\ref{th:main}, but extend the strategy of proof, shown here in its simplest version. Though we do not claim that our result for constant equilibria is significantly new, we have not found it in the literature. In particular, the classical result of \cite[Chapter~4]{Li} is not proved under the sharp spectral assumption, but instead uses a sufficient condition, designed to be able to close the argument by a direct energy estimate. However, we stress that even for the present sharp result one may reasonably argue that a nicer proof could be obtained by replacing spectral theory and Green functions arguments with Fourier analysis. Our present technical choice is purely motivated by the versatility of the designed strategy, having in mind its extension to the proof of Theorem~\ref{th:main}.

To state the result, let us consider $\bU_0\in\R^n$ a zero of $g$, $g(\bU_0)=0_{\,\R^n}$ and set  $\A:=D_UA(\bU_0)\in \mathscr{M}_n(\R)$ and $\G:=D_Ug(\bU_0)\in \mathscr{M}_n(\R)$. We assume that \eqref{edp_intro} is strictly hyperbolic near $\bU_0$, that is, we assume that $\A$ has $n$ distinct real eigenvalues. 

\begin{thm}\label{th:constant}
Assume that there exists $\alpha_0>0$ such that the spectrum of the operator $-\A\partial_x+\G$ (acting on $L^2$ with maximal domain) is contained in 
\[
\{\,\lambda\in\C\,;\,\Re(\lambda)<-\alpha_0\,\}\,.
\]
For any $0<\alpha<\alpha_0$, there exist positive $C_0$ and $\epsilon_0$ such that for any $V_0\in H^2(\R;\R^n)$ with $\|V_0\|_{H^2(\R)}\leq \epsilon_0$, there exist $V\in\mathscr{C}^0(\R_+;H^2(\R;\R^n))\cap \mathscr{C}^1(\R_+;H^1(\R;\R^n))$ such that 
\[
(t,x)\mapsto \bU_0+V(t,x)
\]
solves \eqref{edp_intro} with initial data $\bU_0+V_0(\cdot)$ and
\begin{align*}
\|V(t,\cdot)\|_{H^2(\R)}
&\leq\,C_0\,e^{-\alpha\,t}\,\|V_0\|_{H^2(\R)}\,,&t\in\R_+\,.
\end{align*}
\end{thm}

In the constant stability problem, there is no loss of generality in assuming that $\A$ is not characteristic, that is, that all eigenvalues of $\A$ are non zero. Indeed, the problem is invariant in its assumptions\footnote{The spectrum does change dramatically but the spectral gap does not.} and conclusions by any change of frame $(t,x)\mapsto (t,x-\sigma\,t)$, $\sigma\in \R$, and the latter replaces $\A$ with $\A-\sigma\,I_n$. In the non characteristic case the domain of $-\A\partial_x+\G$ is simply $H^1$. 

From now on, we make the non-characteristic assumption.

\subsection{High-frequency analysis}\label{s:hf-const}

The final argument combines two types of estimates,
\begin{enumerate}
\item estimates on the linearized evolution on one hand, applied on a Duhamel formulation,
\item nonlinear high-frequency damping estimates on the other hand, applied on the original formulation. 
\end{enumerate}
Obviously the former are very directly related to the spectral stability assumption but they are insufficient to conclude by themselves, in the present quasilinear context. In turn, the latter, used to complete the former, are not readily connected to the spectral stability assumption and our first task is precisely to obtain the relevant pieces of information. These are derived from an inspection of the spectrum in the high-frequency regimes. This turns out to be also useful to derive bounds on the linearized evolution, that, in one form or the other, require a \emph{uniform} control of resolvent operators.

To prepare concrete asymptotic expansions, we introduce an invertible diagonal matrix $\D:=\mathrm{diag}(d_j)\in \mathscr{M}_n(\R)$ and an invertible matrix $P\in \mathscr{M}_n(\R)$ such that $\mathbb{A}=P^{-1}\D P$, with $d_1<\cdots<d_j<\cdots<d_n$. To motivate the asymptotic analysis we gather, in advance, the elements used in the nonlinear estimates. 

\begin{lem}\label{l:hf-spec}
Under the assumptions of Theorem~\ref{th:constant},
\begin{enumerate}
\item there exists $\cQ \in \mathscr{M}_n(\R)$ such that $\D^{-1}P\G P^{-1}-\left[ \D^{-1},\cQ \right]=\D^{-1}\Gamma$, where $\Gamma=\mathrm{diag}(\gamma_j) \in \mathscr{M}_n(\R)$ is the diagonal part of $P\G P^{-1}$;  
\item for $j=1,\cdots,n$, there holds $\gamma_j\leq-\alpha_0$.
\end{enumerate}
\end{lem}
In the foregoing lemma, $[A,B]=AB-BA$ stands for the commutator of two matrices $A,B\in \mathscr{M}_n(\R)$. The existence of such a $\cQ$ follows readily from the fact that all the eigenvalues of $\D^{-1}$ are distinct. The content of the lemma is the upper bound on $\gamma_j$.

Instead of rushing at the proof of the lemma, we show how it arises from spectral asymptotics. Thus, for $\lambda\in\C$ with $\Re(\lambda)\geq -\alpha_0$ and $F\in L^2(\R;\C^n)$  we seek for an expansion of solutions to 
\bqq
\lambda \widetilde{V}+ \A \partial_x\widetilde{V}=\G \widetilde{V}+F\,.
\label{edplap}
\eqq
in the regime when the spectral parameter $\lambda$ goes to infinity, $|\lambda|\to\infty$. We can equivalently rewrite \eqref{edplap} as
\bqs
\lambda(P \widetilde{V})+ \D \partial_x(P\widetilde{V})=P\G P^{-1} (P\widetilde{V})+PF,
\eqs
which also reads
\bqs
\partial_x(P\widetilde{V})=-\lambda\D^{-1}(P \widetilde{V}) + \D^{-1}P\G P^{-1} (P\widetilde{V})+\D^{-1}PF.
\eqs
Now, we may use the fact that eigenvalues of $\D$ are distinct to diagonalize the latter equation at a higher order with respect to $\lambda$. This may be carried out at an arbitrary order; for related concrete computations, see \cite{BNRZ,JNRZ2,BJNRZ2}.

To begin with, we introduce $Q_\lambda:=I_n+\dfrac{1}{\lambda}\cQ$ with $\cQ\in \mathscr{M}_n(\R)$ to be fixed later. When $|\lambda|$ is sufficiently large, the matrix $Q_\lambda$ is invertible and \eqref{edplap} takes the form
\bqs
\partial_x(Q_\lambda^{-1}P\widetilde{V})=-\lambda Q_\lambda^{-1} \D^{-1}Q_\lambda(Q_\lambda^{-1}P \widetilde{V}) + Q_\lambda^{-1}\D^{-1}P\G P^{-1}Q_\lambda ( Q_\lambda^{-1} P\widetilde{V})+Q_\lambda^{-1} \D^{-1}Q_\lambda\left(Q_\lambda^{-1}PF\right).
\eqs
Next, we observe that 
\bqs
\lambda \D^{-1}Q_\lambda = Q_\lambda(\lambda \D^{-1})+\left[\lambda \D^{-1},Q_\lambda \right] = Q_\lambda(\lambda \D^{-1})+\left[ \D^{-1},\cQ \right]\,.
\eqs
As a consequence, when $\lambda$ is large, upon introducing the new unknown $W:=Q_\lambda^{-1}P\widetilde{V}$ we arrive at
\bqs
\partial_xW=\left(-\lambda  \D^{-1} +\D^{-1}P\G P^{-1}-\left[ \D^{-1},\cQ \right]+\frac{1}{\lambda} N(\lambda)   \right) W + Q_\lambda^{-1} \D^{-1}Q_\lambda\left(Q_\lambda^{-1}PF\right),
\eqs
for some matrix $N(\lambda)\in \mathscr{M}_n(\R)$ uniformly bounded with respect to $\lambda$ large. At this stage, it is natural to choose $\cQ$ as in Lemma~\ref{l:hf-spec}, using that the eigenvalues of $\D$ are distinct.  

We then arrive at
\bqq
\partial_xW=\MM_\lambda W + Q_\lambda^{-1} \D^{-1}Q_\lambda\left(Q_\lambda^{-1}P F\right)\,,
\label{odeW}
\eqq
with 
\bqs
\MM_\lambda:=-\lambda  \D^{-1} +\D^{-1}\Gamma+\frac{1}{\lambda} N(\lambda)\,,
\eqs
$\Gamma$ diagonal. Anticipating on the sign of $\gamma_j$, stated in Lemma~\ref{l:hf-spec} but still to be proved, we set $\rho_j:=-\gamma_j$. The eigenvalues $\mu_j(\lambda)$ of $\M_\lambda$ expand as 
\begin{align*}
\mu_j(\lambda)&\stackrel{|\lambda|\rightarrow +\infty}{=}\mu_j^\infty(\lambda)+ \mathcal{O}\left(\frac{1}{|\lambda|}\right)\,,&
\mu_j^\infty(\lambda)&:=- \frac{\lambda}{d_j}-\frac{\rho_j}{d_j}\,,
\end{align*}
and have corresponding spectral projectors  $\Pi_j(\lambda)$ expanding as 
\bqs
\Pi_j(\lambda)
\stackrel{|\lambda|\rightarrow +\infty}{=}
\Pi_j^0 + \mathcal{O}\left(\frac{1}{|\lambda|}\right)\,,
\eqs
where $\Pi_j^0$ is simply the standard projection onto the $j$th vector of the canonical basis.

In particular, for each $j$, for $\xi\in\R$ sufficiently large, there exists $\lambda_{j,\xi}$ such that $\mu_j(\lambda_{j,\xi})=\mbi\,\xi$, and 
\[
\lambda_{j,\xi}\,\stackrel{|\xi|\rightarrow +\infty}{=}\,-\rho_j-\mbi\,d_j\xi
+\mathcal{O}\left(\frac{1}{|\xi|}\right)\,.
\]
Now, to deduce the upper bound on $\gamma_j$, we only need to check that the condition $\Im(\mu_j(\lambda))=0$ for some $j$ implies that $\lambda$ belongs to the spectrum of $-\A\partial_x+\G$. The latter claims stems from the fact that this provides a solution $\widetilde{V}$ to \eqref{edplap} with $F\equiv 0_{\R^n}$, of trigonometric monomial type $e^{\mbi\,\xi\,\cdot}$, whose cut-off approximation yields a sequence $(\widetilde{V}_k)_{k\in\N}$ valued in $H^1$ such that
\[
\frac{\|(\lambda+\A\partial_x-\G)\widetilde{V}_k\|_{L^2}}{\|\widetilde{V}_k\|_{L^2}}
\stackrel{k\to+\infty}{\longrightarrow} 0\,.
\]
See details of a related computation in Lemma~2 of the Appendix to \cite[Chapter~5]{Henry}, or Proposition~2.1 in \cite[Section~2.1]{DR2}. Let us stress that this final part of the argument is classical, and does not use the asymptotic expansion. We shall apply it repeatedly without mention from now on, especially to bound contributions from intermediate spectral frequencies $\lambda$.

This achieves the proof of Lemma~\ref{l:hf-spec}. Moreover, for any $\alpha<\alpha_0$, there exists $M>0$ such that if $\lambda$ is such that $\Re(\lambda)\geq -\alpha$ and $|\lambda|\geq M$, each $\Re(\mu_j(\lambda))$ has the sign of $-d_j$. As a consequence, introducing 
\begin{align*}
\mathcal{J}_s&=\left\{ j \in \llbracket 1,n\rrbracket ~|~ d_j>0 \right\}\,,&
\mathcal{J}_u&=\left\{ j \in \llbracket 1,n\rrbracket ~|~ d_j<0 \right\}\,,
\end{align*}
for such a $\lambda$ and $F\in L^2$, one can represent the $H^1$ solution to \eqref{odeW} as
\begin{align*}
W(x)&= \sum_{j\in \mathcal{J}_s} \int_{-\infty}^x e^{\mu_j(\lambda)(x-y)}\Pi_j(\lambda) Q_\lambda^{-1} \D^{-1}Q_\lambda\left(Q_\lambda^{-1}P F(y)\right)\md y\\
&~~~ - \sum_{j\in \mathcal{J}_u} \int_x^{+\infty} e^{\mu_j(\lambda)(x-y)}\Pi_j(\lambda) Q_\lambda^{-1} \D^{-1}Q_\lambda\left(Q_\lambda^{-1}P F(y)\right)\md y \,.
\end{align*}
Rephrased differently, for such a $\lambda$, we have obtained the following Green kernel representation of the resolvent operator
\bqs
\widetilde{V}(x)=\int_\R \K_\lambda(x-y)F(y)\md y\,,
\eqs
with 
\bqs
\K_\lambda(x)=\left\{\begin{array}{lc}
\sum_{j\in \mathcal{J}_s} e^{\mu_j(\lambda)x}\left(P^{-1}\Pi_j^0\,\D^{-1}P+  \mathcal{O}\left(\frac{1}{|\lambda|}\right)\right), & x>0, \\
- \sum_{j\in \mathcal{J}_u} e^{\mu_j(\lambda)x}\left(P^{-1}\Pi_j^0\,\D^{-1}P+  \mathcal{O}\left(\frac{1}{|\lambda|}\right)\right), & x<0\,.\\
\end{array}
\right.
\eqs

In what follows, we use the outcome of the present subsection to obtain, on one hand, bounds on the linearized dynamics, on the other hand, nonlinear energy estimates.

\subsection{Linear stability}\label{s:lin-const}

To begin with, we go on with Green kernel studies so as to prove the following linear asymptotic stability.

\begin{prop}\label{pr:stab-lin-contant}
Under the assumptions of Theorem~\ref{th:constant}, for any $0<\alpha<\alpha_0$, there exists $C>0$ such that for any $V_0\in L^2(\R)$, there exists a unique solution $V\in\mathscr{C}^0(\R_+;L^2(\R))$ to
\bqq
\partial_t V+ \A \partial_xV=\G V,
\label{edplin}
\eqq
with $V(0,\cdot)=V_0$, and, moreover,
\begin{align*}
\|V(t,\cdot)\|_{L^2(\R)}
&\leq\,C\,e^{-\alpha\,t}\,\|V_0\|_{L^2(\R)}\,,&t\in\R_+\,.
\end{align*}
\end{prop}

Since System~\eqref{edplin} is constant-coefficient, thus commutes with derivatives, the foregoing proposition is readily transferred into an $H^k$ result.

We stress again that the proof we give for Proposition~\ref{pr:stab-lin-contant} is by no means the shortest one, but is motivated by further extensions. Indeed, in the $L^2$ context, it is expendient to use an isometry type result, either through Fourier representation, or through the Gearhart-Pr\"uss theorem. The latter requires a uniform bound on the resolvent operator. On any compact set of the spectral plane, the bound follows from a continuity argument. The bound outside some compact set may be derived through Young's inequality from the above Green kernel representation, which yields a uniform bound on $\| \K_\lambda \|_{L^1(\R)}$. Alternatively, one may obtain such a uniform bound through an energy estimate at the spectral level, similar to nonlinear estimates detailed below.

However, an extension to the $L^p$-setting or to the Riemann-shock stability problem of the Hilbert-type arguments would be very cumbersome, if possible at all. 

We study solutions to \eqref{edplin} with data $V_0$ through the Green kernel representation 
\[
V(t,x)=\langle\K^t(x-\cdot);V_0\rangle\,,
\]
where $\langle\,\cdot\,;\,\cdot\,\rangle$ denotes the duality bracket, and the time-evolution Green kernel is obtained from spectral Green kernels through
\begin{equation}\label{eq:inverse-Laplace}
\K^t(\cdot)
\,=\,\frac{1}{2\mbi\pi}\int_{\eta-\mbi\infty}^{\eta+\mbi\infty} e^{\lambda t} \K_\lambda(\cdot) \md \lambda
\end{equation}
where $\eta$ is arbitrary in $(-\alpha_0,+\infty)$. The foregoing integral is an improper integral valued in distributions on $\R$. We mention that, sometimes, instead of using a duality bracket, we shall abuse rigorous notation and write 
\[
V(t,x)=\int_\R\K^t(x-y)\,V_0(y)\md y\,.
\]
The desired $L^2\to L^2$ estimate for the linearized evolution would stem from Young's inequality if one could show that $(\lambda,x)\mapsto \K_\lambda(x)$ belongs to $L^1(\alpha+\mbi\R;L^1(\R))$. Yet this can not happen and, indeed, $\K^t$ is not an $L^1$-function. Instead, we identify explicitly the most singular parts of $\K^t$ and use the above crude argument to bound the reminder part. To reach this stage, one needs an expansion of $\K_\lambda$ up to order $\lambda^{-2}$.

To provide such an expansion, we introduce notation for next correctors of spatial spectral elements
\begin{align*}
\mu_j(\lambda)&\stackrel{|\lambda|\to +\infty}{=}\mu_j^\infty(\lambda)+\frac{\mu_j^1}{\lambda}+\mathcal{O}\left(\frac{1}{|\lambda|^2}\right)\,,\\
\Pi_j(\lambda)&\stackrel{|\lambda|\to +\infty}{=}
\Pi_j^0 + \frac{1}{\lambda}\Pi_j^1+ \mathcal{O}\left(\frac{1}{|\lambda|^2}\right)\,,
\end{align*}
and leading-order part of spectral Green kernels
\bqs
\K_\lambda^\infty(x):=\left\{\begin{array}{lc}
\sum_{j\in \mathcal{J}_s} e^{\mu_j^\infty(\lambda)x}P^{-1}\Pi_j^0\D^{-1}P, & x>0, \\
- \sum_{j\in \mathcal{J}_u} e^{\mu_j^\infty(\lambda)x}P^{-1}\Pi_j^0\D^{-1}P, & x<0.\\
\end{array}
\right.
\eqs
Inserting higher expansions, one derives
\bqs
\K_\lambda(x)=\K_\lambda^\infty(x)+\frac{1}{\lambda}\K_\lambda^1(x)+\frac{x}{\lambda}\K_\lambda^{1,1}(x)+\frac{1}{\lambda^2}\K_\lambda^2(x),
\eqs
where
\bqs
\K_\lambda^1(x):=\left\{\begin{array}{lc}
\sum_{j\in \mathcal{J}_s} e^{\mu_j^\infty(\lambda)x}P^{-1}\left(\Pi^1_j + [\cQ,\Pi^0_j] \right)\D^{-1}P, & x>0, \\
- \sum_{j\in \mathcal{J}_u} e^{\mu_j^\infty(\lambda)x}P^{-1}\left(\Pi^1_j + [\cQ,\Pi^0_j] \right)\D^{-1}P, & x<0,\\
\end{array}
\right.
\eqs
and
\bqs
\K_\lambda^{1,1}(x):=\left\{\begin{array}{lc}
\sum_{j\in \mathcal{J}_s} \mu_j^1 e^{\mu_j^\infty(\lambda)x}P^{-1}\Pi_j^0\D^{-1}P, & x>0, \\
- \sum_{j\in \mathcal{J}_u} \mu_j^1 e^{\mu_j^\infty(\lambda)x}P^{-1}\Pi_j^0\D^{-1}P, & x<0\,.\\
\end{array}
\right.
\eqs

To identify the contributions of the most singular part, we need to partly compute \eqref{eq:inverse-Laplace}. To do so, we point out that \eqref{eq:inverse-Laplace} is simply a way to inverse the Laplace transform, $\lambda\mapsto \K_\lambda(\cdot)$ being the Laplace transform of $t\mapsto \K^t(\cdot)$, and, indeed, we only need to recognize some classical Laplace transforms. However, when doing so, it is useful to switch between the point of view, useful to carry out linear estimates, seeing the temporal Green kernel as the continuous map $t\mapsto \K^t(\cdot)$ valued in distributions on $\R$, and the point of view, also available here thanks to the non-charateristic assumption and useful in Laplace inversions, seeing it as the continuous map $x\mapsto \K^{\cdot}(x)$ valued in distributions on $\R_+$. A typical such identification is that, for $d\neq0$, on the former hand one considers $t\mapsto d\,\boldsymbol{\delta}_{d\,t}(\cdot)$, whereas on the latter hand one manipulates $x\mapsto \delta_{\frac{x}{d}}(\cdot)$, with convential notation that $\delta_{t_0}$, $t_0\in\R_+$, and $\boldsymbol{\delta}_{x_0}$, $x_0\in\R$, denote Dirac masses respectively at $t_0$ and $x_0$. We use different symbols for respective Dirac masses, precisely to prevent any confusion when switching from one point of view to the other. With this in mind, the following identities stem immediately from classical knowledge of Laplace transforms of Dirac masses and indicator functions. For any $t>0$, $1\leq j\leq n$, $\eta>-\rho_j$,
\begin{align*}
\frac{1}{2\mbi\pi}\int_{\eta-\mbi\infty}^{\eta+\mbi\infty} e^{\lambda t}
e^{-(\lambda+\rho_j)\frac{\cdot}{d_j}}\,\chi_{\R_+}\left(\frac{\cdot}{d_j}\right)\,\md \lambda
&=e^{-\rho_j\,t}\,d_j\,\boldsymbol{\delta}_{d_j\,t}(\cdot)\,,\\
\frac{1}{2\mbi\pi}\int_{\eta-\mbi\infty}^{\eta+\mbi\infty} e^{\lambda t}
e^{-(\lambda+\rho_j)\frac{\cdot}{d_j}}\,\chi_{\R_+}\left(\frac{\cdot}{d_j}\right)\,
\frac{\md \lambda}{\lambda+\rho_j}
&=e^{-\rho_j\,\frac{\cdot}{d_j}}\,\chi_{[0,t]}\left(\frac{\cdot}{d_j}\right)\,,\\
\end{align*}
where $\chi_A$ denotes the indicator function of the set $A$. For the reader unfamiliar with Laplace transforms, we stress that in order to check the foregoing claims, it is sufficient to invoke the one-to-one character of the Laplace transform and compute Laplace transforms of right-hand terms, a straightforward task. 

Now, to prove Proposition~\ref{pr:stab-lin-contant}, we fix $\eta=-\alpha$, $\alpha\in(0,\alpha_0)$, and, for some sufficiently large $R$, we split $\K^t(\cdot)$ as 
\begin{align*}
\K^t(\cdot)
&\,=\,
\frac{1}{2\mbi\pi}\int_{\eta-\mbi R}^{\eta+\mbi R} e^{\lambda t} \K_\lambda(\cdot) \md \lambda
+\frac{1}{2\mbi\pi}\int_{\eta+\mbi(\R\setminus [-R,R])} e^{\lambda t}
\frac{1}{\lambda^2}\K_\lambda^2(\cdot)\md \lambda\\
&\quad
+\frac{1}{2\mbi\pi}\int_{\eta+\mbi(\R\setminus [-R,R])}
e^{\lambda t} \left(\K_\lambda(\cdot)
-\frac{1}{\lambda^2}\K_\lambda^2(\cdot)\right)\md \lambda\,.
\end{align*}
Contributions from the first line are bounded directly in $L^1(\R)$ as sketched above. The remaining task is to check that contributions from the second line fit the above explicit computations up to reminders directly bounded in $L^1(\R)$. This follows from the facts that in integrals over $\eta+\mbi(\R\setminus [-R,R])$ one may replace $\lambda^{-1}$ with $(\lambda+\rho_j)^{-1}$ (for the relevant $\rho_j$) up to an $\mathcal{O}(\lambda^{-2})$ reminder, and that once this is done, one may complete integrals over $\eta+\mbi(\R\setminus [-R,R])$ in integrals over $\eta+\mbi\R$. The outcome is as follows.

\begin{lem}\label{l:Green-HF-constant}
Under the assumptions of Theorem~\ref{th:constant}, for any $0<\alpha<\alpha_0$, there exists $C>0$ such that
\begin{align*}
\|\K^t(\cdot)-\K^t_\infty(\cdot)\|_{L^1(\R)}
&\leq\,C\,e^{-\alpha\,t}\,,&t\in\R_+\,,
\end{align*}
where 
\[
\K^t_\infty(\cdot) =  \sum_{j=1}^n d_j\,\boldsymbol{\delta}_{d_j\,t}(\cdot)\ e^{-\rho_j t}\,P^{-1}\Pi_j^0\D^{-1}P\,.
\]
\end{lem} 

This is sufficient to conclude the proof of Proposition~\ref{pr:stab-lin-contant}. We recall that what we have proved also yields with the same argument $L^p$ bounds. In the reverse direction, we point out that even when one chooses to use Fourier-type arguments instead of spectral arguments, when $L^p$ bounds, $p\neq2$, are needed, a similar ammount of work is needed and it is classical to use a decomposition similar to Lemma~\ref{l:Green-HF-constant}. See for instance high-frequency estimates in \cite{Hoff-Zumbrun_1,Hoff-Zumbrun_2,Rodrigues-compressible}.

\subsection{Nonlinear stability}

To prove Theorem~\ref{th:constant}, we need to bound solutions $V$ to 
\begin{equation}\label{eq:constant}
\partial_t V+D_UA(\underline{U}_0+V)\partial_x V=g(\underline{U}_0+V)
\end{equation}
starting from $V_0\in H^2(\R)$ sufficiently small at $t=0$. Classical local well-posedness theory provides a solution $V\in \mathscr{C}^0([0,T_*);H^2(\R;\R^n))\cap \mathscr{C}^1([0,T_*);H^1(\R;\R^n))$, for some maximal existence time $T_*\in(0,+\infty]$, with continuous dependence on $V_0$, and blow-up criterion expressed in terms of $W^{1,\infty}$ topology. System~\eqref{eq:constant} is equivalently written as
\begin{align*}
\partial_t V+(\A\partial_x-\G)V
&=-\partial_x(A(\underline{U}_0+V)-A(\underline{U}_0)-D_UA(\underline{U}_0)V)\\
&\quad+g(\underline{U}_0+V)-g(\underline{U}_0)-D_Ug(\underline{U}_0)V\,.
\end{align*}

As a consequence, applying Duhamel formula and Proposition~\ref{pr:stab-lin-contant} one deduces that, under the assumptions of Theorem~\ref{th:constant}, for any $\alpha<\alpha_0$, there exists a constant $C$ such that if $t<T_*$ is such that $\max_{s\in[0,t]}\|V(s,\cdot)\|_{L^\infty}\leq 1$ then
\bqs
\| V(t,\cdot) \|_{L^2}\leq Ce^{-\alpha t} \| V_0 \|_{L^2}
+C\int_0^t e^{-\alpha (t-s)} 
\|V(s,\cdot)\|_{W^{1,\infty}} \| V(s,\cdot)\|_{L^2}\md s\,.
\eqs

Since the $H^2$ norm controls the $W^{1,\infty}$ norm, to achieve the proof it is sufficient to prove the following $H^2$ high-frequency damping estimate.

\begin{prop}\label{pr:hf-damping-constant}
Under the assumptions of Theorem~\ref{th:constant}, for any $0<\alpha'<\alpha_0$, there exists $C>0$ and $\epsilon>0$ such that for any $H^2$ maximal solution $V$ to \eqref{eq:constant} and any $t$ such that $\max_{s\in[0,t]}\|V(s,\cdot)\|_{W^{1,\infty}}\leq \epsilon$, there holds 
\bqs
\| V(t,\cdot) \|_{H^2}^2\leq Ce^{-2\alpha' t} \| V_0 \|_{H^2}^2
+C\int_0^t e^{-2\alpha' (t-s)}\| V(s,\cdot)\|_{L^2}^2\md s\,.
\eqs
\end{prop}

Indeed, choosing $\alpha'>\alpha$ and combining with the above $L^2$ bound and a Gr\"onwall argument, this shows that as long as the $W^{1,\infty}$ remains small, the $H^2$ norm is exponentially damped with rate $\alpha$ and the $W^{1,\infty}$ norm is kept even smaller. From this, a continuity argument derives that the latter conclusions hold globally in time.

The rest of the section is devoted to the proof of Proposition~\ref{pr:hf-damping-constant}. Here we use crucially the content of Lemma~\ref{l:hf-spec} and mimick the classical damping argument originating in the seminal work of Kawashima \cite{Kawashima}, the latter being a specific instance of the more recently formalized family of hypocoercive damping estimates. For a short introduction to the classical Kawashima theory and further references to the corresponding extensive literature, we refer the reader to \cite[Appendix~A]{R_HDR}. Alternatively, for a similar purpose the reader may also consult \cite{Nguyen,CrinBarat}.

To present first a streamlined version of the computation, we begin by proving a linearized $H^1$ damping estimate, that is, a bound 
\bqs
\| V(t,\cdot) \|_{H^1}^2\leq Ce^{-2\alpha' t} \| V_0 \|_{H^1}^2
+C\int_0^t e^{-2\alpha' (t-s)}\| V(s,\cdot)\|_{L^2}^2\md s\,,
\eqs
for $V$ solving $\partial_t V+(\A\partial_x-\G)V=0_{\R^n}$, $V(0,\cdot)=V_0$. This is done by introducing an $H^1$ functional exponentially dissipated up to $L^2$ remainders. Consider 
\bqs
\mathscr{E}_{lin}(V):=\frac{1}{2}\left\|P\partial_x V \right\|_{L^2(\R)}^2+\langle \bQ PV,P\partial_x V\rangle_{L^2(\R)}+\frac{\vartheta}{2}\left\| P V \right\|_{L^2(\R)}^2,
\eqs
where $\bQ\in\mathscr{M}_n(\R)$ is fixed explicitly below and $\vartheta>0$ is taken sufficiently large to ensure that $\sqrt{\mathscr{E}_{lin}(\cdot)}$ is equivalent to $\|\cdot\|_{H^1(\R)}$. Since $P\A P^{-1}$ is diagonal, when $V$ solves the linearized equation, an integration by parts yields
\begin{align*}
\frac{1}{2}\frac{\md}{\md t}\|PV\|_{L^2}^2
&=\langle P \G P^{-1} P V,PV\rangle_{L^2}\,,&
\frac{1}{2}\frac{\md}{\md t}\|P\partial_xV\|_{L^2}^2
&=\langle P \G P^{-1} P\partial_xV,P\partial_xV\rangle_{L^2}\,,
\end{align*}
whereas
\begin{align*}
\frac{\md}{\md t}\langle \bQ PV,P\partial_x V\rangle_{L^2} &=\langle[\D,\bQ]\,P\partial_x V,P\partial_x V\rangle_{L^2}+\langle \bQ P\G P^{-1} PV,P\partial_x V\rangle_{L^2}
+\langle \bQ PV,P\G P^{-1} P\partial_x V\rangle_{L^2}\,.
\end{align*}
At this stage, we recall from Lemma~\ref{l:hf-spec} that there exists $\cQ \in \mathscr{M}_n(\R)$ such that $P\G P^{-1}-(\cQ-\D\,\cQ\,\D^{-1})=\Gamma$, with $\Gamma=\mathrm{diag}(\gamma_j) \in \mathscr{M}_n(\R)$, where for $j=1,\cdots,n$, $\gamma_j\leq-\alpha_0$. Setting  $\bQ:=\cQ\D^{-1}$ so that $P\G P^{-1}+[\D,\bQ]=\Gamma$, one deduces that when $V$ solves the linearized equation,
\begin{align*}
\frac{\md}{\md t}\mathscr{E}_{lin}(V)
-\alpha_0\,\left\|P\partial_x V \right\|_{L^2(\R)}^2
&\lesssim \|PV\|_{L^2}\,\|P\partial_xV\|_{L^2}+\vartheta \|PV\|_{L^2}^2\,,
\end{align*}
which implies, through Young inequality,
\begin{align*}
\frac{\md}{\md t}\mathscr{E}_{lin}(V)
-\alpha'\,\left\|P\partial_x V \right\|_{L^2(\R)}^2
&\lesssim (1+\vartheta)\|PV\|_{L^2}^2\,,
\end{align*}
thus also
\begin{align*}
\frac{\md}{\md t}\mathscr{E}_{lin}(V)
-2\,\alpha'\,\mathscr{E}_{lin}(V)
&\lesssim (1+\vartheta)\|V\|_{L^2}^2\,.
\end{align*}
Integrating and using the equivalence of $\sqrt{\mathscr{E}_{lin}(\cdot)}$ with the standard $H^1$ norm achieves the proof of the claim.

To extend the $H^1$ damping estimate to the nonlinear system, we introduce $\P(\cdot)$ a smooth map, defined on a neighborhood of $\bU_0$, and valued in invertible matrices, such that $\P(\bU_0)=P$ and, for any $U$, $\P(U)\,D_UA(U)\,\P(U)^{-1}$ is diagonal. Consider
\bqs
\mathscr{E}_1(V):=\frac{1}{2}\left\| \P(\bU_0+V) \partial_x V \right\|_{L^2(\R)}^2+\langle \bQ PV,P\partial_x V\rangle_{L^2(\R)}+\frac{\vartheta}{2}\left\| P V \right\|_{L^2(\R)}^2,
\eqs
with $\bQ$ as above and $\vartheta$ possibly larger. Arguing essentially as in the linear case, one derives that if $V$ solves \eqref{eq:constant}, then, as long as the $L^\infty$ norm of $V$ is kept sufficiently small, $\mathscr{E}_1(V)$ is equivalent to $\|V\|_{H^1}^2$ and 
\begin{align*}
\frac{\md}{\md t}\mathscr{E}_1(V)
-\alpha_0\,\left\|P\partial_x V \right\|_{L^2(\R)}^2
&\lesssim 
\|V\|_{W^{1,\infty}}\,\|P\partial_xV\|_{L^2}^2\\
&\qquad+\left(1+\|V\|_{W^{1,\infty}}\right)
\left(\|PV\|_{L^2}\,\|P\partial_xV\|_{L^2}+\|PV\|_{L^2}^2\right)\,,
\end{align*}
(with an implicit dependence on $\vartheta$), from which an $H^1$ damping estimate with rate $\alpha'$ is obtained as long as the $W^{1,\infty}$ norm of $V$ is kept sufficiently small.

With $H^1$ bounds in hands, in order to close the proof of Proposition~\ref{pr:hf-damping-constant} thus of Theorem~\ref{th:constant}, it is sufficient to carry out similar computations for the functional 
\bqs
\mathscr{E}_2(V):=\frac{1}{2}\left\| \P(\bU_0+V) \partial_x^2 V \right\|_{L^2(\R)}^2+\langle \bQ P\partial_xV,P\partial_x^2 V\rangle_{L^2(\R)}+\frac{\vartheta}{2}\left\| P  \partial_xV \right\|_{L^2(\R)}^2\,.
\eqs

Note that a smallness constraint on the $W^{1,\infty}$ norm arises already when proving an $H^1$ damping estimate but the constraint does not get stringer when going to $H^k$ damping estimates, $k\geq 1$. 

\begin{rmk}
As already mentioned, our proof may readily be adapted to yield a $W^{2,p}$-stability result, $1\leq p<\infty$, and with a little more work it also provides a $W^{s,p}$-stability result, when $s>1+1/p$, $1<p<\infty$. When $2<p<\infty$, this relaxes the localization constraint on initial perturbations. To give a clue on required changes, we point out that when deriving a linearized $W^{1,p}$ damping, $1<p<\infty$, $\mathscr{E}_{lin}(V)$ should be replaced with  
\bqs
\frac{1}{p}\sum_j\int_\R|(P\partial_x V)_j|^p\md x
+\sum_j\int_\R (\bQ PV)_j |(P\partial_x V)_j|^{p-1} \sgn((P\partial_x V)_j)
\md x+\frac{\vartheta}{p}\sum_j\int_\R|(PV)_j|^p\md x\,.
\eqs
\end{rmk}

\section{Initial boundary value problems}\label{s:IBVP}

In this section, we introduce one last preparatory problem before tackling the proof of Theorem~\ref{th:main}.

We consider the following system of equations
\bqq
\left\{
\begin{split}
\partial_t U+ \partial_x(A(U))&=g(U),\qquad \textrm{on }\R_+\times\R_+,\\
B[U(\cdot,0)]&=\varphi,\qquad\quad\textrm{on }\R_+,\\
U(0,\cdot)&=U_0,\qquad\quad\textrm{on }\R_+\,,
\end{split}
\right.
\label{edpx>0}
\eqq
with $A$ and $g$ as in previous sections and $B$ a smooth map from $\R^n$ to $\R^p$ for some integer $p$.

We pick $\bU_0\in\R^n$ such that $(t,x)\mapsto \bU_0$ is a stationnary solution to \eqref{edpx>0}, that is, such that $g(\bU_0)=0_{\R^n}$ and $B[\bU_0]=0_{\R^p}$. In vague terms, our goal is to prove that if such a solution is spectrally stable with a spectral gap, then, when $U_0$ is sufficiently close to $\bU_0$, $\varphi$ is sufficiently small in a space encoding exponential time decay, and $U_0$ and $\varphi$ are sufficiently compatible with each other, solutions to \eqref{edpx>0} decay exponentially in time to $\bU_0$. To be more precise, we prove such a result when, besides spectral stability, we assume that near $\bU_0$ the system is strictly hyperbolic and the boundary $x=0$ is non characteristic, that is, when moreover $\A:=D_UA(\bU_0)$ have $n$ distinct real nonzero eigenvalues.

Before stating the corresponding result, since we expect less readers to be familiar with spectrum for IBVPs, let us be more explicit on the latter. Set $\G:=D_U g(\bU_0)$ and $\B:=D_U B(\bU_0)$. For any spectral parameter $\lambda\in\C$, we say that $\lambda$ does not belong to the spectrum of \eqref{edpx>0} linearized about $(t,x)\mapsto \bU_0$ provided that for any $F\in L^2(\R_+;\C^n)$ and any $F_0\in\C^p$, there exists a unique $\widetilde{V}\in H^1(\R_+;\C^n)$ solving
\bqq
\left\{
\begin{split}
\left(\lambda I_n+\A\partial_x-\G\right)\widetilde{V}&=F \qquad 
\textrm{on }\R_+,\\
\B\widetilde{V}(0)&=F_0.
\end{split}
\right.
\label{resx>0}
\eqq
Note that with such a definition, if the spectrum is not the entire complex plane, $\B$ is onto (and in particular $p\leq n$).

Then, as in the situation of Theorem~\ref{th:main}, we observe that prescribing some regularity structure on solutions is essentially equivalent to imposing compatibility constraints on data. Thus, as in the introduction, to keep the main statement of the present section as streamlined as possible, we introduce beforehand the terminology that perturbation $(V_0,\varphi)$ are $H^2$-compatible if $B[\bU_0+V_0(0)]=\varphi(0)$ and
\begin{align*}
D_U\,B(\bU_0+V_0)\left(-\partial_x(A(\bU_0+V_0))+g(\bU_0+V_0)\right)&=\varphi'(0)\,.
\end{align*}

\begin{thm}\label{th:IBVP}
Assume that there exists $\alpha_0>0$ such that the spectrum of the linearization of \eqref{edpx>0} about $\bU_0$ is contained in 
\[
\{\,\lambda\in\C\,;\,\Re(\lambda)<-\alpha_0\,\}\,.
\]
For any $0<\alpha<\alpha_0$, there exist positive $C_0$ and $\epsilon_0$ such that for any $(V_0,\varphi)\in H^2(\R_+;\R^n)\times BUC^2(\R_+;\R^p)$ with $\|V_0\|_{H^2(\R_+)}+\|e^{\alpha\,\cdot}\varphi\|_{W^{2,\infty}(\R_+)}\leq \epsilon_0$, there exist $V\in\mathscr{C}^0(\R_+;H^2(\R_+;\R^n))\cap \mathscr{C}^1(\R_+;H^1(\R_+;\R^n))$ such that 
\[
(t,x)\mapsto \bU_0+V(t,x)
\]
solves \eqref{edpx>0} with initial data $\bU_0+V_0(\cdot)$ and boundary data $\varphi$ and
\begin{align*}
\|V(t,\cdot)\|_{H^2(\R_+)}
&\leq\,C_0\,e^{-\alpha\,t}\,\left(\|V_0\|_{H^2(\R_+)}
+\|e^{\alpha\,\cdot}\varphi\|_{W^{2,\infty}(\R_+)}\right)\,,&t\in\R_+\,.
\end{align*}
\end{thm}

In the foregoing theorem, we use notation $BUC^k(\Omega)$, with $k\in\N$ and $\Omega$ connected, to denote the set of functions whose derivatives up to order $k$ are bounded and uniformly continuous. 

\begin{rmk}
Assuming exponential decay of $\varphi$ is obviously necessary to prove exponential decay of $U$ since $B[U(t,0)]=\varphi(t)$. Yet, when a weaker decay is assumed, say $(\varphi,\varphi',\varphi'')(t)\lesssim 1/\omega(t)$ for some decay rate function $\omega$, the proof also yields nonlinear stability, this time with decay rate $\omega$, provided that for some $0<\alpha<\alpha_0$, $1\lesssim \omega(t)\lesssim e^{\alpha\,t}$, 
\begin{align*}
\int_0^t e^{-\alpha\,(t-s)}\,\,\frac{\md\,s}{\omega(s)^2}
&\lesssim \frac{1}{\omega(t)}\,,&
\textrm{and}&&
\int_0^t e^{-2\,\alpha\,(t-s)}\,\,\frac{\md\,s}{\omega(s)^2}
&\lesssim \frac{1}{\omega(t)^2}\,.
\end{align*}
This holds for instance for $\omega$ given by $\omega(t)=(1+t)^\epsilon$, for some $\epsilon\geq0$.
\end{rmk}

\begin{rmk}
For the present IBVP, our proof may also be extended to yield an $L^p$-based result, and to higher-regularity results. Note however that short-time persistence of higher regularity requires more compatibility conditions. A simple way to ensure compatibility at any order is to assume that $\varphi$ is supported away from the initial time and $V_0$ is supported away from the spatial origin.
\end{rmk}

\subsection{Preliminary spectral analysis}

The conclusions of Lemma~\ref{l:hf-spec} hold as they are under the assumptions of Theorem~\ref{th:IBVP} and the arguments expounded to prove it provide a few more useful conclusions for the situation at hand.

To begin with, arguing as in Section~\ref{s:hf-const}, one recovers the classical fact that exponential dichotomy is necessary to invertibility of the spectral problem: for any $\lambda$ not in the spectrum, the matrix $\LL(\lambda):=\A^{-1}\left(\G-\lambda I_n\right)$ possesses no purely imaginary eigenvalues. When this holds, denoting by $\Pi_{s}(\lambda)$ and $\Pi_{u}(\lambda)$ the spectral projections associated with stable and unstable spaces of $\LL(\lambda)$, one then derives readily that invertibility of the spectral problem is equivalent to invertibility of $\B\arrowvert_{\Ran(\Pi_s(\lambda))}$. Moreover, then, solutions to \eqref{resx>0} are obtained through the matrix-valued Green kernel representation
\bqs
\widetilde{V}(x)=\K_\lambda^{\bc}(x)F_0
+\int_{\R^+} \K_\lambda^{\hom}(x,y)F(y)\md y\,,
\eqs
with
\begin{align*}
\K_\lambda^{\bc}(x)&=e^{\LL(\lambda)x} \left(\B\arrowvert_{\Ran(\Pi_s(\lambda))}\right)^{-1}\,,&
\K_\lambda^{\hom}(x,y)&=\K_\lambda^{\dir}(x-y)+\K_\lambda^{\rfl}(x,y)\,,&
\end{align*}
where
\begin{align*}
\K_\lambda^{\dir}(x)&=\left\{
\begin{array}{lc}
e^{\LL(\lambda)x}\,\Pi_s(\lambda)\,\A^{-1},& 0<x\,,\\
- e^{\LL(\lambda)x}\,\Pi_u(\lambda)\,\A^{-1},&x<0\,,
\end{array}
\right.\\
\K_\lambda^{\rfl}(x,y)&=e^{\LL(\lambda)x} \left(\B\arrowvert_{\Ran(\Pi_s(\lambda))}\right)^{-1}\,\B e^{-\LL(\lambda)y}\,\Pi_u(\lambda)\,\A^{-1}\,.
\end{align*}
The foregoing exponents stand respectively for \emph{boundary condition}, \emph{homogeneous}, \emph{direct} and \emph{reflected}. The derivation of the latter formula is essentially immediate from
\begin{align*}
\widetilde{V}(x)= e^{\LL(\lambda)x}\Pi_s(\lambda)\widetilde{V}(0)
+\int_0^x e^{\LL(\lambda)(x-y)}\,\Pi_s(\lambda)\,\A^{-1}\,F(y)\md y
-\int_x^{+\infty} e^{\LL(\lambda)(x-y)}\,\Pi_u(\lambda)\,\A^{-1}\,F(y) \md y\,,
\end{align*}
and spectral boundary condition $\B \Pi_s(\lambda)\widetilde{V}(0)=F_0-\B \Pi_u(\lambda)\widetilde{V}(0)$.

Inserting the high-frequency expansions of Section~\ref{s:hf-const} in the foregoing representation provides directly the high-frequency expansions required for the present analysis. Explicitly the leading-order part of $\K_\lambda^{\dir}(x)$ is given by
\[
\left\{\begin{array}{lc}
\sum_{j\in \mathcal{J}_s} e^{\mu_j^\infty(\lambda)x}P^{-1}\,\Pi_j^0\,\D^{-1}P, & x>0, \\
- \sum_{j\in \mathcal{J}_u} e^{\mu_j^\infty(\lambda)x}P^{-1}\,\Pi_j^0\,\D^{-1}P, & x<0,\\
\end{array}
\right.
\]
the one of $\K_\lambda^{\bc}(x)$ is
\[
\sum_{j\in \mathcal{J}_s} e^{\mu_j^\infty(\lambda)x}P^{-1}\,\Pi_j^0\,P\,
\left(\B\arrowvert_{\Ran(\Pi_s^\infty)}\right)^{-1}
\]
and the one of $\K_\lambda^{\rfl}(x,y)$ is
\[
\sum_{j\in \mathcal{J}_s}\sum_{\ell\in \mathcal{J}_u}  e^{\mu_j^\infty(\lambda)x}\,e^{-\mu_\ell^\infty(\lambda)y}\,P^{-1}\,\Pi_j^0\,P\,
\left(\B\arrowvert_{\Ran(\Pi_s^\infty)}\right)^{-1}\B\,P^{-1}\,\Pi_\ell^0\,\D^{-1}P
\]
where $\Pi_s^\infty:=\sum_{j\in \mathcal{J}_s}P^{-1}\Pi_j^0P$ is the projection on incoming characteristics of the linearized system. The next order of the expansions is likewise available and it is indeed also required to apply our arguments providing linear stability.

Note that, as expected, the Lopatinski\u{\i} condition that $\B\arrowvert_{\Ran(\Pi_s^\infty)}$ is invertible arises as a consequence of the spectral stability assumption (and not as an extra preliminary assumption).

\subsection{Linear stability}

As in Section~\ref{s:constant}, we first prove a linear asymptotic stability result.

\begin{prop}\label{pr:stab-lin-IBVP}
Under the assumptions of Theorem~\ref{th:IBVP}, for any $0<\alpha<\alpha'<\alpha_0$, there exist positive $C$ and $C'$ such that for any $V_0\in L^2(\R_+)$ and $\varphi\in BUC^0(\R_+)$, there exists a unique solution $V\in\mathscr{C}^0(\R_+;L^2(\R_+))$ to
\bqq
\left\{
\begin{split}
\partial_t V+\A \partial_xV&=\G V, \qquad\quad\textrm{on }\R_+\times\R_+,\\
\B V(\cdot,0)&=\varphi,\qquad\quad\textrm{on }\R_+,\\
V(0,\cdot)&=V_0,\qquad\quad\textrm{on }\R_+\,,
\end{split}
\right.
\label{linx>0}
\eqq
and, moreover,
\begin{align*}
\|V(t,\cdot)\|_{L^2(\R_+)}
&\leq\,
C\,e^{-\alpha\,t}\,\|V_0\|_{L^2(\R_+)}
+C'\,\int_0^t\,e^{-\alpha'\,(t-s)}\,\|\varphi(s)\|\,\md s\\
&\leq\,C\,e^{-\alpha\,t}\,
\,\left(\|V_0\|_{L^2(\R_+)}
+\|e^{\alpha\,\cdot}\varphi\|_{L^{\infty}(\R_+)}\right)\,,&t\in\R_+\,.
\end{align*}
\end{prop}
For the sake of clarity, we recall the classical observation that in such low regularity results the existence of traces is not derived from classical trace theorems for Sobolev spaces but arises as a consequence of the evolution equation and its non characterisc character.

We derive Proposition~\ref{pr:stab-lin-IBVP} from Green kernel representation
\begin{align*}
V(t,x)&=\langle\K^x_{\bc}(t-\cdot);\varphi\rangle+\langle\K^t_{\hom}(x,\cdot);V_0\rangle\\
&=\langle\K^x_{\bc}(t-\cdot);\varphi\rangle
+\langle\K^t_{\rfl}(x,\cdot);V_0\rangle
+\langle\K^t_{\dir}(x-\cdot);V_0\rangle\,,
\end{align*}
where time-dependent Green kernels are obtained from their spectral counterparts as in Section~\ref{s:lin-const}.

Proposition~\ref{pr:stab-lin-IBVP} follows from the following lemma.

\begin{lem}\label{l:Green-HF-IBVP}
Under the assumptions of Theorem~\ref{th:IBVP}, for any $0<\alpha<\alpha_0$, there exists $C>0$ such that
\begin{align*}
\|\K^t_{\dir}(\cdot)-\K^t_{\dir,\infty}(\cdot)\|_{L^1(\R)}
+\|\K^{\,\cdot\,}_{\bc}(t)-\K^{\,\cdot\,}_{\bc,\infty}(t)\|_{L^1\cap L^\infty(\R_+)}\\
+\sup_x \|\K^t_{\rfl}(x,\cdot)-\K^t_{\rfl,\infty}(x,\cdot)\|_{L^1(\R_+)}\\
+\sup_y \|\K^t_{\rfl}(\cdot,y)-\K^t_{\rfl,\infty}(\cdot,y)\|_{L^1(\R_+)}
&\leq\,C\,e^{-\alpha\,t}\,,&t\in\R_+\,,
\end{align*}
with
\begin{align*}
\K^t_{\dir,\infty}(\cdot)&=  \sum_{j=1}^n\chi_{\R_+^*}(d_j\,t)\ d_j\,\boldsymbol{\delta}_{d_j\,t}(\cdot)\ e^{-\rho_j t}\,P^{-1}\Pi_j^0\D^{-1}P\,,\qquad t>0\,,\\
\K^{x}_{\bc,\infty}(\cdot)
&=\sum_{j\in \mathcal{J}_s}\,\delta_{\frac{x}{d_j}}(\cdot)\ e^{-\rho_j\,\frac{x}{d_j}}\,P^{-1}\Pi_j^0P\,\left(\B\arrowvert_{\Ran(\Pi_s^\infty)}\right)^{-1}\,,\qquad x>0\,,
\end{align*}
and
\begin{align*}
\K^t_{\rfl,\infty}(x,\cdot)
&=\sum_{j\in \mathcal{J}_s}\sum_{\ell\in \mathcal{J}_u} 
\Bigg[\chi_{\R_+^*}\left(t-\frac{x}{d_j}\right)
\,|d_\ell|\,\boldsymbol{\delta}_{|d_\ell|\left(t-\frac{x}{d_j}\right)}(\cdot)\ 
e^{-\rho_j\,\frac{x}{d_j}}\,e^{-\rho_\ell\,\left(t-\frac{x}{d_j}\right)}\\
&\qquad\qquad\qquad
\times\,P^{-1}\,\Pi_j^0\,P\,\left(\B\arrowvert_{\Ran(\Pi_s^\infty)}\right)^{-1}\B\,P^{-1}\,\Pi_\ell^0\,\D^{-1}P\Bigg]\,,\qquad t>0\,,\ x>0\,.
\end{align*}
\end{lem} 

We omit the proof of Lemma~\ref{l:Green-HF-IBVP} as essentially identical to the one of Lemma~\ref{l:Green-HF-constant}, but give some details on how to deduce from it Proposition~\ref{pr:stab-lin-IBVP}. The contribution of $\K^t_{\dir}$ is estimated as the whole dynamics is in Proposition~\ref{pr:stab-lin-contant}. The contribution of $\K^t_{\rfl}$ is estimated by noticing that an $L^1\to L^1$ bound stems from the $L^\infty_yL^1_x$ control, whereas an\footnote{Actually a $BUC^0\to BUC^0$ bound.} $L^\infty\to L^\infty$ bound stems from the $L^\infty_xL^1_y$ control, hence the $L^2\to L^2$ bound by interpolation. Concerning the contribution of $\K^{\,\cdot\,}_{\bc}(t)$, we simply point out that it is sufficient to apply Lemma~\ref{l:Green-HF-IBVP} with decay rate $\alpha'$.

\subsection{Nonlinear stability}

As in the proof of Theorem~\ref{th:constant}, our proof of Theorem~\ref{th:IBVP} is concluded by a continuity argument on maximal solutions given by the standard local well-posedness theory. On the latter, besides classical references \cite{Li-Yu,Bressan,Metivier,Benzoni-Serre}, we refer the reader to the recent \cite{Audiard} for optimal regularity results and a concise introduction.

One half of the required estimates is directly given by applying Proposition~\ref{pr:stab-lin-IBVP} to
\begin{align*}
\partial_t V+(\A\partial_x-\G)V
&=-\partial_x(A(\underline{U}_0+V)-A(\underline{U}_0)-D_UA(\underline{U}_0)V)\\
&\quad+g(\underline{U}_0+V)-g(\underline{U}_0)-D_Ug(\underline{U}_0)V\,,\\
\B V(\cdot,0)
&=\,\varphi-\left(B[\underline{U}_0+V(\cdot,0)]
-B(\underline{U}_0)
-D_U B(\underline{U}_0)V\right)\,,\\
V(0,\cdot)&=V_0\,.
\end{align*}
This yields that under the assumptions of Theorem~\ref{th:IBVP}, for any $\alpha<\alpha_0$, there exists a constant $C$ such that if the solution is defined on $[0,t]$ and satisfies $\max_{s\in[0,t]}\|V(s,\cdot)\|_{L^\infty}\leq 1$ then
\bqs
\| V(t,\cdot) \|_{L^2}\leq Ce^{-\alpha t}
\,\left(\|V_0\|_{L^2}
+\|e^{\alpha\,\cdot}\varphi\|_{L^{\infty}}\right)
+C\int_0^t e^{-\alpha (t-s)} 
\|V(s,\cdot)\|_{W^{1,\infty}} \| V(s,\cdot)\|_{L^2}\md s\,.
\eqs
To prove the latter claim, the only new ingredient we have used besides Proposition~\ref{pr:stab-lin-IBVP} is
\[
\|V(s,0)\|\lesssim \|V(s,\cdot)\|_{W^{1,\infty}}^{1/2} \| V(s,\cdot)\|_{L^2}^{1/2}\,,
\]
that follows from Sobolev embedding inequality 
$\|V(s,\cdot)\|_{L^\infty}\lesssim \|\nabla V(s,\cdot)\|_{L^{\infty}}^{1/3} \| V(s,\cdot)\|_{L^2}^{2/3}$ (and the trivial bound $\|V(s,0)\|\leq \|V(s,\cdot)\|_{L^\infty}$).

To achieve the proof of Theorem~\ref{th:IBVP}, it is then sufficient to prove the following $H^2$ high-frequency damping estimate.

\begin{prop}\label{pr:hf-damping-IBVP}
Under the assumptions of Theorem~\ref{th:IBVP}, for any $0<\alpha'<\alpha_0$, there exists $C>0$ and $\epsilon>0$ such that for any $H^2$ maximal solution $V$ to \eqref{edpx>0} and any $t$ such that $\max_{s\in[0,t]}\|V(s,\cdot)\|_{W^{1,\infty}}\leq \epsilon$, there holds 
\bqs
\| V(t,\cdot) \|_{H^2}^2\leq Ce^{-2\alpha' t}\| V_0 \|_{H^2}^2
+C\int_0^t e^{-2\alpha' (t-s)}\,
\left(\|(\varphi,\varphi',\varphi'')(s)\|^2
+\| V(s,\cdot)\|_{L^2}^2\right)\md s\,.
\eqs
\end{prop}

As in the proof of Proposition~\ref{pr:hf-damping-constant}, the core of the proof of Proposition~\ref{pr:hf-damping-IBVP} is already present in the derivation of a linearized $H^1$ damping estimate, that is, a bound 
\bqs
\| V(t,\cdot) \|_{H^1}^2\leq Ce^{-2\alpha' t} \| V_0 \|_{H^1}^2
+C\int_0^t e^{-2\alpha' (t-s)}
\left(\|(\varphi,\varphi')(s)\|^2+\|V(s,\cdot)\|_{L^2}^2\right)\md s\,.
\eqs
for $V$ solving $\partial_t V+(\A\partial_x-\G)V=0_{\R^n}$, $V(0,\cdot)=V_0$, $\B V(\cdot,0)=\varphi$. For the sake of exposition simplicity, we begin by proving such a linearized estimate. By a density-continuity argument, one may recover the general $H^1$ case from the subcase when $(V_0,\varphi)$ are smooth and $H^2$-compatible, or even from the subcase when $(V_0,\varphi)$ are smooth and compatible at any order. We thus focus on the latter.

The key difference with Proposition~\ref{pr:hf-damping-constant} is that energy estimates involve boundary terms. Roughly speaking, as far as high-frequency damping estimates are concerned, outgoing characteristics are associated with dissipative boundary terms thus help in closing estimates, whereas incoming characteristics yield boundary terms to be controlled by outgoing boundary terms through the Lopatinski\u{\i} condition. This is already seen on $L^2$ estimates. Indeed, with $\Pi_u^\infty:=I-\Pi_s^\infty$, when $V$ solves the announced linearized problem,
\begin{align*}
\frac{1}{2}\frac{\md}{\md t}\left\|P\,\Pi_{s}^\infty V\right\|^2_{L^2(\R_+)}(t)
&-\frac{1}{2}\langle \D P\Pi_{s}^\infty V(t,0), P\Pi_{s}^\infty V(t,0)\rangle
\lesssim \left\|V(t,\cdot)\right\|^2_{L^2(\R_+)}\,,\\
\frac{1}{2}\frac{\md}{\md t}\left\|P\,\Pi_{u}^\infty V\right\|^2_{L^2(\R_+)}(t)
&-\frac{1}{2}\langle \D P\Pi_{u}^\infty V(t,0), P\Pi_{u}^\infty V(t,0)\rangle
\lesssim \left\|V(t,\cdot)\right\|^2_{L^2(\R_+)}\,,\\
\end{align*}
whereas $-P^*\D P$ is positive definite on $\Ran(\Pi_u^\infty)$ since 
\[
-(\Pi_u^\infty)^*P^*\D P\Pi_u^\infty\,=\,
-\sum_{j\in \mathcal{J}_u} P^*\Pi_j^0\D\Pi_j^0P
\geq (\min_{j\in \mathcal{J}_u}|d_j|)\,(\Pi_u^\infty)^*\Pi_u^\infty\,.
\]
Since $\B\arrowvert_{\Ran(\Pi_s^\infty)}$ is invertible and $\B(\Pi_s^\infty V(t,0))=\varphi(t)-\B(\Pi_u^\infty V(t,0))$, one deduces that for some $c>0$ and any $\theta'$ sufficiently large
\begin{align*}
\frac{1}{2}\frac{\md}{\md t}
\left(\theta'\left\|P\,\Pi_{u}^\infty V\right\|^2_{L^2}
+\left\|P\,\Pi_{s}^\infty V\right\|^2_{L^2}\right)(t)
&+c\,\theta'\|V(t,0)\|^2
\lesssim (1+\theta')\,
\left(\|\varphi(t)\|^2+\left\|V(t,\cdot)\right\|^2_{L^2(\R_+)}\right)\,.
\end{align*}
The foregoing estimate encodes that when $\theta'$ is sufficiently large, $\theta'\,(P\,\Pi_{u}^\infty)^*\,P\,\Pi_{u}^\infty+(P\,\Pi_{s}^\infty)^*\,P\,\Pi_{s}^\infty$ is a symmetrizer for which boundary conditions are dissipative. The presence of an $L^2$-norm in the right-hand side is due to the fact that $\G$ possesses no particular structure for this symmetrizer.

To incorporate similar elements in the estimate of $\partial_xV$, we need to identify corresponding boundary conditions. This is achieved by differentiating with respect to the time variable the boundary equation to obtain
\[
\B \A \partial_xV(t,0)\,=\,\B\G V(t,0)\,-\,\varphi'(t)\,.
\]
Note that $\A$ is invertible on $\Ran(\Pi_s^\infty)$ so that the foregoing computations also yield that for some $c>0$ and any $\theta'$ sufficiently large
\begin{align*}
\frac{1}{2}&\frac{\md}{\md t}
\left(\theta'\left\|P\,\Pi_{u}^\infty\partial_xV\right\|^2_{L^2}
+\left\|P\,\Pi_{s}^\infty\partial_xV\right\|^2_{L^2}\right)(t)\\
&+c\,\theta'\|\partial_xV(t,0)\|^2
-\theta'\langle P\Pi_{u}^\infty\partial_x V(t,\cdot), P\Pi_{u}^\infty\G \partial_x V(t,\cdot)\rangle_{L^2}
-\langle P\Pi_{s}^\infty \partial_xV(t,\cdot), P\Pi_{s}^\infty\G \partial_xV(t,\cdot)\rangle_{L^2}\\
&\lesssim (1+\theta')\,\left(\|(\varphi,\varphi')(t)\|^2+\,\left\|V(t,0)\right\|^2\right)\,.
\end{align*}

Using that $\A$ commutes with $\Pi_{s}^\infty$ and $\Pi_{u}^\infty$, or, equivalently, that $\D$ commutes with $P\Pi_{s}^\infty P^{-1}$ and $P\Pi_{u}^\infty P^{-1}$, the Kawashima compensator part of the estimates used in the proof of  Proposition~\ref{pr:hf-damping-constant} may also be split according to outgoing and incoming characteristics. Indeed, when $V$ solves the linearized problem under study, for $\#\in\{s,u\}$,
\begin{align*}
\frac{\md}{\md t}\langle P\Pi_{\#}^\infty P^{-1} \bQ PV,P\Pi_{\#}^\infty P^{-1} P\partial_x V\rangle_{L^2}(t)&-\langle P\Pi_{\#}^\infty P^{-1} [\D,\bQ]\,P\partial_x V(t,\cdot),P\Pi_{\#}^\infty P^{-1}\,P\partial_x V(t,\cdot)\rangle_{L^2}\\
&\lesssim 
 \|V(t,0)\|\,\|\partial_x V(t,0)\|
+\|V(t,\cdot)\|_{L^2}\|\partial_x V(t,\cdot)\|_{L^2}\,.
\end{align*}
Thus, choosing again $\bQ$ such that $P\G P^{-1}+[\D,\bQ]=\Gamma$ and using that $\Gamma$ commutes with $P\Pi_{s}^\infty P^{-1}$ and $P\Pi_{u}^\infty P^{-1}$, and $-\Gamma\geq \alpha_0\,I_n$, one deduces for
\begin{align*}
\mathscr{E}_{lin}(V)
&:=
\theta'\left(\frac{1}{2}\left\|P\,\Pi_{u}^\infty\partial_x V \right\|_{L^2}^2+\langle P\,\Pi_{u}^\infty P^{-1}\bQ PV,P\,\Pi_{u}^\infty\partial_x V\rangle_{L^2}+\frac{\vartheta}{2}\left\|P\,\Pi_{u}^\infty V \right\|_{L^2}^2\right)\\
&\quad 
+\frac{1}{2}\left\|P\,\Pi_{s}^\infty\partial_x V \right\|_{L^2}^2+\langle P\,\Pi_{s}^\infty P^{-1}\bQ PV,P\,\Pi_{s}^\infty\partial_x V\rangle_{L^2}+\frac{\vartheta}{2}\left\|P\,\Pi_{s}^\infty V \right\|_{L^2}^2
\end{align*}
that 
\begin{align*}
\frac{\md}{\md t}\mathscr{E}_{lin}(V)(t)
&-\alpha_0\,\left(\theta'\,\left\|P\,\Pi_{u}^\infty\partial_x V \right\|_{L^2}^2
+\left\|P\,\Pi_{s}^\infty\partial_x V \right\|_{L^2}^2\right)\\
&\lesssim_{\theta,\theta'}\,\|(\varphi,\varphi')(t)\|^2
+\|V(t,\cdot)\|_{L^2}^2+\|V(t,\cdot)\|_{L^2}\|\partial_x V(t,\cdot)\|_{L^2}\,,
\end{align*}
provided that $\theta$ and $\theta'$ are sufficiently large. From here, the proof of the linearized estimate is achieved as is the corresponding bound expounded along the proof of Proposition~\ref{pr:hf-damping-constant}.

The extension to the nonlinear problem follows also the strategy carried out to prove Proposition~\ref{pr:hf-damping-constant}, mainly replacing $P$, $\Pi_{u}^\infty$, $\Pi_{s}^\infty$ associated with $\A=D_UA(\bU_0)$ with nonlinear versions associated with $D_UA(\bU_0+ V)$, when $V$ is small. We skip corresponding details.

\section{Stability of Riemann shocks}

We finally turn our attention to the stability of Riemann shocks so as to prove Theorem~\ref{th:main}. The main difference with the initial boundary value problem is that the position of the boundary (at the shock) is free. This results in the introduction of a phase shift $\psi$ tracking the shock position. 

\subsection{Preliminary spectral analysis}

We begin with spectral considerations. To emphasize similarities with the initial boundary value problem of the previous section, we introduce a map $\B$ from $\R\times \R^n\times\R^n$ to $\R^n$ defined by
\[
\B(\Phi,W_+,W_-)\,=\,-\Phi\,[\bU]_0+\A_+\,W_+-\A_-W_-
\]
where $\A_-:=D_U A(\bU_-)-\sigma\,I_n$ and $\A_+:=D_U A(\bU_+)-\sigma\,I_n$. Conversely, to untangle the triplet, we shall use canonical projections $\bI_0$, $\bI_+$ and $\bI_-$, defined by $\bI_0 (\Phi,W_+,W_-)=\Phi$, $\bI_+ (\Phi,W_+,W_-)=W_+$ and $\bI_- (\Phi,W_+,W_-)=W_-$, and associated canonical sections $\bI^0$, $\bI^+$ and $\bI^-$, defined by $\bI^0\Phi=(\Phi,0_{\,\R^n},0_{\,\R^n})$, $\bI^+W_+=(0,W_+,0_{\,\R^n})$ and $\bI^-W_-=(0,0_{\,\R^n},W_-)$.

A large body of the spectral analysis of the foregoing sections is also directly applicable to constant-coefficient operators associated respectively with $\bU_+$ and $\bU_-$, including exponential dichotomies and high-frequency expansions. We shall denote with subscripts ${}_-$ and ${}_+$ the corresponding objects: $\A_-$, $\A_+$, $\B_-$, $\B_+$, $\Gamma_-$, $\Gamma_+$, $\bQ_-$, $\bQ_+$, $\Pi_{s,-}^\infty$, $\Pi_{s,+}^\infty$,...

In particular, for any $\lambda$ not in the spectrum, the matrices $\LL_+(\lambda)=\A_+^{-1}\left(\G_+-\lambda I_n\right)$ and $\LL_-(\lambda)=\A_-^{-1}\left(\G_--\lambda I_n\right)$  possess no purely imaginary eigenvalues. Moreover, when this holds and $\lambda\neq0$, denoting by $\Pi_{s,\pm}(\lambda)$ and $\Pi_{u,\pm}(\lambda)$ the corresponding spectral projections associated with stable and unstable spaces of $\LL_{\pm}(\lambda)$, one then derives readily that invertibility of the spectral problem is equivalent to invertibility of $\B\arrowvert_{\C\times\Ran(\Pi_{s,+}(\lambda))\times\Ran(\Pi_{u,-}(\lambda))}$. Moreover, then, solutions to 
\bqq
\left\{
\begin{split}
\left(\lambda I_n+\A_\pm \partial_x-\G_\pm\right)\widetilde{V}_\pm&=F_\pm, \qquad 
\textrm{on }\R_\pm\,,\\
\quad\B(\lambda \widetilde{\psi},\widetilde{V}_+(0),\widetilde{V}_-(0))&=F_0\,,
\end{split}
\right.
\label{linspectral}
\eqq
are obtained through the matrix-valued Green kernel representation
\begin{align*}
\widetilde{V}_+(x)&=\K_\lambda^{\bc_+}(x)F_0
+\int_{\R^+}\left(\K_\lambda^{\dir_+}(x-y)+\K_\lambda^{\rfl_{++}}(x,y)\right)F_+(y)\md y
+\int_{\R^-}\K_\lambda^{\rfl_{+-}}(x,y)F_-(y)\md y\,,\\
\widetilde{V}_-(x)&=\K_\lambda^{\bc_-}(x)F_0
+\int_{\R^-}\left(\K_\lambda^{\dir_-}(x-y)+\K_\lambda^{\rfl_{--}}(x,y)\right)F_-(y)\md y
+\int_{\R^+}\K_\lambda^{\rfl_{-+}}(x,y)F_+(y)\md y\,,\\
\lambda\widetilde{\psi}&=\K_\lambda^{\bc_0}F_0
+\int_{\R^+}\K_\lambda^{\trp_{+}}(y)F_+(y)\md y
+\int_{\R^-}\K_\lambda^{\trp_{-}}(y)F_-(y)\md y\,,
\end{align*}
with
\begin{align*}
\K_\lambda^{\bc_{\pm}}(x)&=e^{\LL_{\pm}(\lambda)x} \bI_{\pm}\,\B^\dagger(\lambda)\,,&
\K_\lambda^{\bc_0}&=\bI_{0}\,\B^\dagger(\lambda)\,,&
\end{align*}
\begin{align*}
\K_\lambda^{\dir_{\pm}}(x)&=\left\{
\begin{array}{lc}
e^{\LL_{\pm}(\lambda)x}\,\Pi_{s,{\pm}}(\lambda)\,\A_{\pm}^{-1},& 0<x\,,\\
- e^{\LL_{\pm}(\lambda)x}\,\Pi_{u,{\pm}}(\lambda)\,\A_{\pm}^{-1},&x<0\,,
\end{array}
\right.
\end{align*}
\begin{align*}
\K_\lambda^{\rfl_{\pm,+}}(x,y)&=e^{\LL_{\pm}(\lambda)x}\bI_{\pm}\,\B^\dagger(\lambda)
\,\B\,\bI^+e^{-\LL_+(\lambda)y}\,\Pi_{u,+}(\lambda)\,\A_+^{-1}\,,\\
\K_\lambda^{\rfl_{\pm,-}}(x,y)&=-e^{\LL_{\pm}(\lambda)x}\bI_{\pm}\,\B^\dagger(\lambda)
\,\B\,\bI^-e^{-\LL_-(\lambda)y}\,\Pi_{s,-}(\lambda)\,\A_-^{-1}\,,
\end{align*}
and
\begin{align*}
\K_\lambda^{\trp_{+}}(y)&=\bI_{0}\,\B^\dagger(\lambda)
\,\B\,\bI^+e^{-\LL_+(\lambda)y}\,\Pi_{u,+}(\lambda)\,\A_+^{-1}\,,\\
\K_\lambda^{\trp_{-}}(y)&=-\bI_{0}\,\B^\dagger(\lambda)
\,\B\,\bI^-e^{-\LL_-(\lambda)y}\,\Pi_{s,-}(\lambda)\,\A_-^{-1}\,,
\end{align*}
where $\B^\dagger(\lambda)$ denotes the inverse of $\B\arrowvert_{\C\times\Ran(\Pi_{s,+}(\lambda))\times\Ran(\Pi_{u,-}(\lambda))}$. We mention that the new exponent stands for \emph{trapped}. 

From the above representations, one readily deduces that the assumption that $0$ is a simple eigenvalue is equivalent to the fact that on one hand exponential dichotomy also holds at $\lambda=0$ and on the other hand $\B\arrowvert_{\C\times\Ran(\Pi_{s,+}(0))\times\Ran(\Pi_{u,-}(0))}$ is invertible. 

In the reverse direction, inserting high-frequency expansions yields that $\B\arrowvert_{\C\times\Ran(\Pi_{s,+}^\infty)\times\Ran(\Pi_{u,-}^\infty)}$ is invertible and we shall denote $\B^\dagger_\infty$ its inverse. Note that, in particular, the number of characteristics incoming into the shock from the left and the right sum to $n-1$. The shock under study is a Lax shock.

\subsection{Linear stability}

We now prove a linear asymptotic stability result.

\begin{prop}\label{pr:stab-lin}
Under the assumptions of Theorem~\ref{th:main}, for any $0<\alpha<\alpha'<\alpha_0$, there exist positive $C$ and $C'$ such that for any $V_0\in L^2\cap BUC^0(\R^*)$, $\psi_0\in\R$ and $\varphi\in BUC^0(\R_+)$ such that there exists $\psi_1\in\R$ such that $\B(\psi_1,V_0(0^+),V_0(0^-))=\varphi(0)$, there exists a unique solution $(V,\psi)\in\mathscr{C}^0(\R_+;L^2\cap BUC^0(\R^*))\times \mathscr{C}^1(\R_+)$ to
\bqq
\left\{
\begin{split}
\partial_t V+\A_+ \partial_xV&=\G_+ V, \qquad\quad\textrm{on }\R_+\times\R_+^*,\\
\partial_t V+\A_- \partial_xV&=\G_- V, \qquad\quad\textrm{on }\R_+\times\R_-^*,\\
\quad\B(\psi',V(\cdot,0^+),V(\cdot,0^-))&=\varphi,\qquad\qquad\ \textrm{on }\R_+,\\
V(0,\cdot)&=V_0,\qquad\qquad\textrm{on }\R^*\,,\\
\psi(0)&=\psi_0\,,
\end{split}
\right.
\label{linsystedp}
\eqq
and, moreover, for any $t\in\R_+$,
\begin{align*}
\|V(t,\cdot)\|_{L^2\cap L^\infty(\R^*)}
&\leq\,
C\,e^{-\alpha\,t}\,\|V_0\|_{L^2\cap L^\infty(\R^*)}
+C'\,\int_0^t\,e^{-\alpha'\,(t-s)}\,\|\varphi(s)\|\,\md s\\
&\leq\,C\,e^{-\alpha\,t}\,
\,\left(\|V_0\|_{L^2\cap L^\infty(\R^*)}
+\|e^{\alpha\,\cdot}\varphi\|_{L^{\infty}(\R_+)}\right)\,,\\
|\psi'(t)|
&\leq\, 
C\,e^{-\alpha\,t}\,\|V_0\|_{L^\infty(\R^*)}
+C'\,\|\varphi(t)\|
+C'\,\int_0^t\,e^{-\alpha'\,(t-s)}\,\|\varphi(s)\|\,\md s\\
&\leq\,C\,e^{-\alpha\,t}\,
\,\left(\|V_0\|_{L^\infty(\R^*)}
+\|e^{\alpha\,\cdot}\varphi\|_{L^{\infty}(\R_+)}\right)\,.
\end{align*}
\end{prop}

In the foregoing proposition, as in \cite{DR1,DR2}, our convention is that when $\Omega$ is not connected, $BUC^k(\Omega)$ denotes the set of functions that are $BUC^k$ on each connected component of $\Omega$. In particular $BUC^0(\R^*)$ cannot be identified with $BUC^0(\R)$.

The proof of Proposition~\ref{pr:stab-lin} follows quite closely the one of Proposition~\ref{pr:stab-lin-IBVP} so that we provide details only about the newest part, the estimate of $\psi'$. Note however that for the first time we do use that our method also provides $L^\infty\to L^\infty$ bounds. Nevertheless we stress that we do so only to state a linearized stability result as satisfactory as possible but a linearized $L^2$ bound on $V$ would be sufficient to close the nonlinear argument and prove Theorem~\ref{th:main}.

From the Green kernel spectral representation, one deduces a time-dependent Green kernel representation
\begin{align*}
\psi'(t)&=\langle\K_{\bc_0}(t-\cdot);\varphi\rangle
+\langle\K^t_{\trp};V_0\rangle\,.
\end{align*}
The $\psi'$-estimate is then deduced from the following lemma, whose proof stems from high-frequency expansions essentially as corresponding lemmas of former sections. With this respect, let us only mention that in the identification of the subprincipal part of $\K_{\bc_0}$, we use that for any $t>0$ and $\eta<0$,
\begin{align*}
\frac{1}{2\mbi\pi}\int_{\eta-\mbi\infty}^{\eta+\mbi\infty} e^{\lambda t}\,
\frac{\md \lambda}{\lambda}
&=0\,.
\end{align*}

\begin{lem}\label{l:Green-HF}
Under the assumptions of Theorem~\ref{th:main}, for any $0<\alpha<\alpha_0$, there exists $C>0$ such that
\begin{align*}
\|\K_{\bc_0}(t)-\K_{\bc_0,\infty}(t)\|
+\|\K^t_{\trp}(\cdot)-\K^t_{\trp,\infty}(\cdot)\|_{L^1\cap L^\infty(\R^*)}
&\leq\,C\,e^{-\alpha\,t}\,,&t\in\R_+\,,
\end{align*}
with
\begin{align*}
\K_{\bc_0,\infty}(\cdot)
&=\delta_0(\cdot)\,\bI_0\,\B^\dagger_\infty\,,
\end{align*}
and
\begin{align*}
\K^t_{\trp,\infty}(\cdot)
&=\sum_{\ell\in \mathcal{J}_{u,+}} 
\,|d_{\ell,+}|\,\boldsymbol{\delta}_{|d_{\ell,+}|t}(\cdot)
\,e^{-\rho_{\ell,+}\,t}\,\bI_0\,\B^\dagger_\infty\,
\B\,\bI^+\,P_+^{-1}\,\Pi_{\ell,+}^0\,\D_+^{-1}P_+\\
&-\sum_{\ell\in \mathcal{J}_{s,-}} 
\,d_{\ell,-}\,\boldsymbol{\delta}_{-d_{\ell,-}t}(\cdot)
\,e^{-\rho_{\ell,-}\,t}\,\bI_0\,\B^\dagger_\infty\,
\B\,\bI^-\,P_-^{-1}\,\Pi_{\ell,-}^0\,\D_-^{-1}P_-\,,\qquad t>0\,.
\end{align*}
\end{lem} 

\subsection{Nonlinear stability}

We finally prove Theorem~\ref{th:main}. Recall that we consider solutions in the form 
\bqs
U(t,x)=\bU(x-\sigma t -\psi(t))+V(t,x-\sigma t -\psi(t)),
\eqs
with $\bU$ and $V(t,\cdot)$ piecewise smooth with discontinuity at $0$, so that they satisfy
\bqs
\left\{
\begin{split}
\partial_t V+ \left(D_UA(\bU+V)-(\sigma +\psi'(t))I_n\right)\partial_x V&=g(\bU+V)\,,
\qquad \textrm{on }\R_+\times\R^*\\
\quad
B(\sigma+\psi',\bU(0^+)+V(\cdot,0^+),\bU(0^-)+V(\cdot,0^-))&=0_{\,\R^n}\,,
\qquad \qquad\ \textrm{on }\R_+
\end{split}
\right.
\eqs
where, to stress similarities with the fixed-boundary problem, we have introduced $B$ defined by
\[
B(\Phi,W_+,W_-)\,=\,-\Phi\,(W_+-W_-)+A(W_+)-A(W_-)\,.
\]

Note that, consistently with the linearized problem studied in the foregoing subsection, there do hold
$B(\sigma,\bU(0^+),\bU(0^-))=0_{\,\R^n}$ and
\begin{align*}
D_UA(\bU)-\sigma\,I_n&\,=\,\chi_{\R_+^*}\,\A_++\chi_{\R_-^*}\,\A_-\,,\\
D_Ug(\bU)&\,=\,\chi_{\R_+^*}\,\G_++\chi_{\R_-^*}\,\G_-\,,\\
D_{(\Phi,W_+,W_-)}B(\sigma,\bU(0^+),\bU(0^-))\,&=\B\,.
\end{align*}

Estimates of the previous section yield that under the assumptions of Theorem~\ref{th:main}, for any $\alpha<\alpha_0$, there exists a constant $C$ such that if the solution is defined on $[0,t]$ and satisfies $\max_{s\in[0,t]}\|V(s,\cdot)\|_{L^\infty}\leq 1$ then
\begin{align*}
\| V(t,\cdot) \|_{L^2\cap L^\infty}
&\leq Ce^{-\alpha t}\,\|V_0\|_{L^2\cap L^\infty}
+C\int_0^t e^{-\alpha (t-s)} 
\left(\|V(s,\cdot)\|_{W^{1,\infty}}+|\psi'(s)|\right)
\left(\| V(s,\cdot)\|_{L^2}+|\psi'(s)|\right)\md s\,,\\
|\psi'(t)|
&\leq Ce^{-\alpha t}\,\|V_0\|_{L^\infty}
+C\left(\|V(t,\cdot)\|_{W^{1,\infty}}+|\psi'(t)|\right)
\left(\| V(t,\cdot)\|_{L^2}+|\psi'(t)|\right)\\
&\qquad+C\int_0^t e^{-\alpha (t-s)} 
\left(\|V(s,\cdot)\|_{W^{1,\infty}}+|\psi'(s)|\right)
\left(\| V(s,\cdot)\|_{L^2}+|\psi'(s)|\right)\md s\,.
\end{align*}
We point out that the estimate on $\psi'$ is quite rough, but sufficient, and that a direct inspection of the specific form of the Rankine-Hugoniot conditions would improve the bound, but in a useless way.

Therefore, to achieve the proof of Theorem~\ref{th:main} by a continuity argument, it is then sufficient to prove the following $H^2$ high-frequency damping estimate.

\begin{prop}\label{pr:hf-damping}
Under the assumptions of Theorem~\ref{th:main}, for any $0<\alpha'<\alpha_0$, there exists $C>0$ and $\epsilon>0$ such that for any $H^2$ maximal solution $(V,\psi)$ to \eqref{edpx>0} and any $t$ such that $\max_{s\in[0,t]}(|\psi'(s)|+|\psi''(s)|+\|V(s,\cdot)\|_{W^{1,\infty}})\leq \epsilon$, there holds 
\bqs
\| V(t,\cdot) \|_{H^2}^2+|\psi''(t)|^2
\leq Ce^{-2\alpha' t}\| V_0 \|_{H^2}^2
+C\int_0^t e^{-2\alpha' (t-s)}\,
\left(\|\psi'(s)\|^2
+\| V(s,\cdot)\|_{L^2\cap L^\infty}^2\right)\md s\,.
\eqs
\end{prop}

The overall strategy and most of technical computations involved in the proof of Proposition~\ref{pr:hf-damping} are identical to the ones expounded along the proof of corresponding propositions of previous sections. Therefore we only provide details about what differ from the latter, that is, about the parts involving the phase position $\psi$. Since some key differences arise at the nonlinear level, we directly discuss the proof of Proposition~\ref{pr:hf-damping}. We recall that when doing so, thanks to a density argument, it is sufficient to consider smoother solutions (arising from more compatible data).

As a preliminary let us stress that when considering smooth waves, the phase shift is not uniquely determined by the dynamics and one may enforce as an extra normalizing condition that the phase shift is low-frequency so that it plays essentially no role in the nonlinear closing in regularity. This is in strong contrast with the present case when the phase shift is uniquely determined and has a limited amount of smoothness.

The main difference with the nonlinear estimate of the previous section is that the nonlinear objects extending $P_+$, $P_-$, $\Pi_{s,\pm}^\infty$, $\Pi_{u,\pm}^\infty$ are designed to be associated with $\A[\Phi,U]:=D_UA(U)-\Phi\,I_n$ and thus depend on both $U$ and $\Phi$, with $U$ to be replaced with $\bU+V(t,\cdot)$ and $\Phi$ to be replaced with $\sigma+\psi'(t)$. Therefore, as far as interior equations are concerned, controlling $(\psi'(t),\psi''(t))$ has essentially the same role as controlling $\|V(t,\cdot)\|_{W^{1,\infty}}$.

To go on with the discussion, let us denote $\Pi_{s}[\Phi,U]$ and $\Pi_{u}[\phi,U]$ the corresponding extensions. Note that when $(\phi,U_+,U_-)$ is sufficiently close to $(\sigma,\bU_+,\bU_-)$, $\B[\Phi,W_+,W_-]:=D_{(\Phi,W_+,W_-)}B(\phi,U_+,U_-)$ restricted to $\R\times \Ran(\Pi_{s}[\phi,U_+])\times \Ran(\Pi_{u}[\phi,U_-])$ is invertible (with a smooth inverse). This provides a control on both incoming characteristics (as in the previous section) and derivatives of $\psi'$, through, at first order,
\begin{align*}
&\B[\sigma+\psi',(\bU+V)(0^+),(\bU_-+V)(0^-)]
\begin{pmatrix}\psi''\\
-\A[\sigma+\psi',(\bU+V)(0^+)]\partial_x V(0^+)\\
-\A[\sigma+\psi',(\bU+V)(0^-)]\partial_x V(0^-)
\end{pmatrix}\\
&\,=\,
-\B[\sigma+\psi',(\bU+V)(0^+),(\bU_-+V)(0^-)]
\begin{pmatrix}0\\
g((\bU+V)(0^+))\\
g((\bU+V)(0^-))
\end{pmatrix}
\end{align*}
(where we have left implicit time dependencies and used column notation to spare some room) and a similar second-order boundary equation.

Up to these points, the proof of Proposition~\ref{pr:hf-damping} is identical to the one of Proposition~\ref{pr:hf-damping-IBVP}. Through a continuity argument this achieves the proof of Theorem~\ref{th:main} except for the part involving $\psi_\infty$. But this one is readily deduced by integration from bounds on $\psi'$, with
\[
\psi_\infty\,=\,\psi_0+\int_0^{\infty}\psi'(s)\,\md s\,.
\]

\appendix

\section{Stability vs. dissipativity}\label{s:example}

In this short section, we elucidate to which extent spectral stability is a strictly larger notion than dissipative symmetrizability, for constant-coefficient hyperbolic systems. 

Without loss of generality we may consider systems $\partial_t V+\A\,\partial_x V=\G\,V$, with $\A$ and $\G$ real, and $\A$ diagonal. Our claim is that for a suitable choice of system, $0$ is exponentially spectrally stable but the system is not strictly dissipatively symmetrizable, in the sense that there does not exist a symmetric positive definite matrix $S$ such that $S\A$ is symmetric and the real part of $S\G$ is negative.

\subsection{$2\times 2$ systems}

We first restrict to systems of two equations
\begin{align*}
\A&=\begin{pmatrix} d_1&0\\0&d_2\end{pmatrix}\,,&
\G&=\begin{pmatrix} a&b\\c&d\end{pmatrix}\,.&
\end{align*}
In this case we show that the two notions do coincide.

Note that if $d_1=d_2$ then $S\A$ is symmetric for any symmetric $S$ so that the result follows from the well-known equivalence for finite-dimensional ODEs. Thus we assume $d_1\neq d_2$.

Let us first use Fourier computations to enforce spectral stability. We want to ensure that there exists $\theta>0$ such that for any $\xi\in\R$, the eigenvalues of $-\mbi\xi\,\A+\G$ have real part less than $-\theta$. By examination of asymptotic expansion in the limit $|\xi|\to\infty$, one readily checks that this is achieved for $\xi$ large provided that $a<0$ and $d<0$. We conclude by examining under which condition no transition can occur when varying $\xi$. To do so, note that for $\tau\in\R$, $\mbi\,\tau$ is an eigenvalue of $-\mbi\xi\,\A+\G$ if and only if
\begin{align*}
(\xi\,d_1+\tau)\,d&=-(\xi\,d_2+\tau)\,a\,,&
(\xi\,d_1+\tau)\,(\xi\,d_2+\tau)\,&=ad-bc\,,
\end{align*}
and, when $ad>0$, this possesses no solution if and only if $ad-bc>0$. Thus exponential spectral stability is equivalent to $a<0$, $d<0$ and $ad-bc>0$.

In turn, one readily checks that the set of symmetric positive definite matrices $S$ such that $S\A$ is symmetric is exactly the set of matrices 
\[
S\,=\,\begin{pmatrix} \alpha_1&0\\0&\alpha_2\end{pmatrix}
\]
with $\alpha_1>0$ and $\alpha_2>0$. For such a matrix $S$, one has
\[
S\G=\begin{pmatrix} \alpha_1\,a&\alpha_1\,b\\
\alpha_2\,c&\alpha_2\,d\end{pmatrix}\,,
\]
whose real part is negative if and only if 
\begin{align*}
\alpha_1\,a&<0&\textrm{and}&&
\alpha_1\,\alpha_2\,a\,d
&>\frac14\,(\alpha_1\,b+\alpha_2\,c)^2\,.
\end{align*}
If $a<0$, $d<0$ and $ad-bc>0$ with $bc\neq0$ then the condition is met with $\alpha_1=|c|$ and $\alpha_2=|b|$. If $a<0$, $d<0$ and $bc=0$, then the condition is met with one of the $\alpha_j$s equal to $1$ ($\alpha_1$ if $b=0$, $\alpha_2$ otherwise) and the other one sufficiently small. 

\subsection{$3\times 3$ systems}

We turn to systems of three equations
\begin{align*}
\A&=\begin{pmatrix} d_1&0&0\\0&d_2&0\\0&0&d_3\end{pmatrix}\,,&
\G&=\begin{pmatrix} a_1&b_3&c_2\\c_3&a_2&b_1\\b_2&c_1&a_3\end{pmatrix}\,.&
\end{align*}
We assume $d_1\neq d_2$, $d_2\neq d_3$ and $d_3\neq d_1$.

High-frequency exponential stability is equivalent to $a_1<0$, $a_2<0$ and $a_3<0$. We make this assumption from now on. Transition at frequency $\xi$ with eigenvalue $\mbi\tau$ is equivalent to
\begin{align*}
a_1\,b_1\,c_1+a_2\,b_2\,c_2+a_3\,b_3\,c_3
&=-(\xi\,d_1+\tau)\,(\xi\,d_2+\tau)\,a_3\\
&\quad-(\xi\,d_1+\tau)\,(\xi\,d_3+\tau)\,a_2\\
&\quad-(\xi\,d_2+\tau)\,(\xi\,d_3+\tau)\,a_1\,,&\\
(\xi\,d_1+\tau)\,(\xi\,d_2+\tau)\,(\xi\,d_3+\tau)\,
&=(\xi\,d_1+\tau)\,(a_2\,a_3-b_1\,c_1)\\
&\quad+(\xi\,d_2+\tau)\,(a_1\,a_3-b_2\,c_2)\\
&\quad+(\xi\,d_3+\tau)\,(a_1\,a_2-b_3\,c_3)\,.
\end{align*}
To break the symmetry, let us assume that $d_3\in[d_1,d_2]$ and set $\theta:=(d_3-d_1)/(d_2-d_1)$. Then the existence of a transition is equivalent to the existence of a pair of real numbers $(X,Y)$ such that
\begin{align*}
a_1\,b_1\,c_1+a_2\,b_2\,c_2+a_3\,b_3\,c_3
&=-XY\,a_3-X\,((1-\theta)X+\theta Y)\,a_2
-Y\,((1-\theta)X+\theta Y)\,a_1\,,&\\
XY\,((1-\theta)X+\theta Y)\,
&=X\,(a_2\,a_3-b_1\,c_1)
+Y\,(a_1\,a_3-b_2\,c_2)
+((1-\theta)X+\theta Y)\,(a_1\,a_2-b_3\,c_3)\,.
\end{align*}
The absence of solution with $X=0$ is equivalent to
\begin{align*}
&&a_1\,b_1\,c_1+a_2\,b_2\,c_2+a_3\,b_3\,c_3&<0\,,\\
\textrm{or}&&
(a_1\,b_1\,c_1+a_2\,b_2\,c_2+a_3\,b_3\,c_3>0
&\quad\textrm{ and }\quad
a_1\,a_3-b_2\,c_2
+\theta\,(a_1\,a_2-b_3\,c_3)\neq0)\,.
\end{align*}
For $X\neq0$, the existence of a $Y$ such that $(X,Y)$ solves the system is equivalent to the existence of a $Z$ such that 
\begin{align*}
X^2\,\left(a_1\,b_1\,c_1+a_2\,b_2\,c_2+a_3\,b_3\,c_3\right)
&=-X^2Z\,a_3-X^2\,((1-\theta)X^2+\theta Z)\,a_2
-Z\,((1-\theta)X^2+\theta Z)\,a_1\,,&\\
Z\,((1-\theta)X^2+\theta Z)\,
&=X^2\,(a_2\,a_3-b_1\,c_1)
+Z\,(a_1\,a_3-b_2\,c_2)
+((1-\theta)X^2+\theta Z)\,(a_1\,a_2-b_3\,c_3)\,.
\end{align*}

Specializing to 
\begin{align*}
d_1&=1\,,&d_2&=3\,,&d_3&=2\,,&\theta&=\frac12\,,&
\end{align*}
\begin{align*}
a_1=a_2=a_3&=-1\,,&
b_2=b_3=c_1=c_2&=1\,,&
b_1=c_3&=-1\,.
\end{align*}
The equations of the system become 
\begin{align*}
Z\,\left(-\frac32\,X^2\,-\,1\right)&=\frac12 X^4\,+2X^2\,
\end{align*}
and
\begin{align*}
\frac12\,Z^2\,
&=Z\,(-\frac12\,X^2+1)+3\,X^2\,.
\end{align*}
This implies
\begin{align*}
\frac12\,\left(\frac12 X^4\,+2X^2\right)^2\,
&=\left(-\frac32\,X^2\,-\,1\right)\,\left(\frac12 X^4\,+2X^2\right)
\,\left(-\frac12\,X^2+1\right)\\
&\quad+3\,X^2\,\left(-\frac32\,X^2\,-\,1\right)^2\,,
\end{align*}
thus
\begin{align*}
\frac14\,X^6+\frac{27}{4} X^4+\frac92 X^2+1\,=\,0\,.
\end{align*}
This is impossible. Therefore one deduces spectral stability for 
\begin{align*}
\A&=\begin{pmatrix} 1&0&0\\0&3&0\\0&0&2\end{pmatrix}\,,&
\G&=\begin{pmatrix} -1&1&1\\-1&-1&-1\\1&1&-1\end{pmatrix}\,.&
\end{align*}

The set of symmetric positive definite matrices $S$ such that $S\A$ is symmetric is exactly the set of matrices 
\[
S\,=\,\begin{pmatrix} \alpha_1&0&0\\0&\alpha_2&0\\0&0&\alpha_3\end{pmatrix}
\]
with $\alpha_1>0$, $\alpha_2>0$ and $\alpha_3>0$. For such a matrix $S$, one has
\[
S\G=\begin{pmatrix} \alpha_1\,a_1&\alpha_1\,b_3&\alpha_1\,c_2\\
\alpha_2\,c_3&\alpha_2\,a_2&\alpha_2\,b_1\\
\alpha_3\,b_2&\alpha_3\,c_1&\alpha_3\,a_3\end{pmatrix}\,,
\]
whose real part is negative if and only if 
\begin{align*}
\alpha_1\,a_1&<0\,,&
\alpha_1\,a_1\,\alpha_3\,a_3
&>\frac14\,(\alpha_1\,c_2+\alpha_3\,b_2)^2\,,&
\end{align*}
and
\begin{align*}
\alpha_1\,a_1\,\alpha_2\,a_2\,\alpha_3\,a_3\,
&+\frac14\,(\alpha_1\,b_3+\alpha_2\,c_3)\,(\alpha_2\,b_1+\alpha_3\,c_1)
\,(\alpha_3\,b_2+\alpha_1\,c_2)\\
&>\frac14\,\alpha_1\,a_1\,(\alpha_2\,b_1+\alpha_3\,c_1)^2
+\frac14\,\alpha_2\,a_2\,(\alpha_3\,b_2+\alpha_1\,c_2)^2
+\frac14\,\alpha_3\,a_3\,(\alpha_1\,b_3+\alpha_2\,c_3)^2\,.
\end{align*}
Note that for our specific above choice, the second condition becomes
\[
0>\frac14\,(\alpha_1-\alpha_3)^2\,,
\] 
which is impossible.

Note moreover that replacing the above choice with
\begin{align*}
\A&=\begin{pmatrix} 1&0&0\\0&3&0\\0&0&2\end{pmatrix}\,,&
\G_\epsilon&=\begin{pmatrix} -1&1&1\\-1&-1&-1\\1+\epsilon&1&-1\end{pmatrix}\,.&
\end{align*}
with $\epsilon>0$ sufficiently small, one even obtains an example for which the zero solution is exponentially stable but for any symmetrizer $S$, the real part of $S\G_\epsilon$ possesses a positive eigenvalue.

\newcommand{\etalchar}[1]{$^{#1}$}


\end{document}